\documentclass[moor,nonblindrev]{informs1} 

\usepackage{lmodern}
\usepackage{xcolor}
\usepackage{amsmath}
\usepackage{amssymb, xfrac,enumerate}

\def\LP{\operatorname{LP}}
\def\IP{\operatorname{IP}}

\global\long\def\ve#1{\boldsymbol{#1}}%
\global\long\def\R{\mathbb{R}}%
\global\long\def\Z{\mathbb{Z}}%
\global\long\def\aa{\ve a}%
\global\long\def\AA{\ve A}%
\global\long\def\BB{\ve B}%
\global\long\def\PP{\ve P}%
\global\long\def\SSv{\ve S}%
\global\long\def\MM{\ve M}%
\global\long\def\TT{\ve T}%
\global\long\def\UU{\ve U}%
\global\long\def\bb{\ve b}%
\global\long\def\cc{\ve c}%
\global\long\def\hh{\ve h}%
\global\long\def\gg{\ve g}%
\global\long\def\rr{\ve r}%
\global\long\def\ee{\mathbf{e}}%
\global\long\def\pp{\ve p}%
\global\long\def\uu{\ve u}%
\global\long\def\vv{\ve v}%
\global\long\def\ww{\ve w}%
\global\long\def\xx{\ve x}%
\global\long\def\Q{\mathbb{Q}}%
\global\long\def\yy{\ve y}%
\global\long\def\zz{\ve z}%
\global\long\def\aalpha{\ve{\alpha}}%
\global\long\def\bbeta{\ve{\beta}}%
\global\long\def\zero{\mathbf{0}}%
\global\long\def\one{\mathbf{1}}%
\global\long\def\vol{\mathrm{vol}}%
\global\long\def\sign{\mathrm{sign}}%
\global\long\def\conv{\mathrm{conv}}%
\global\long\def\ccone{\mathcal{C}}%
\global\long\def\poly{\mathcal{P}}%
\global\long\def\norm#1{\left\Vert #1\right\Vert }%
\global\long\def\spindle{\mathcal{S}}%
\global\long\def\qoly{\mathcal{Q}}%
\global\long\def\prox{\kappa}%
\global\long\def\objfn{\aalpha}%

\global\long\def\AAhat{\hat{\AA\;}\negthickspace}%

\global\long\def\rot{\tau}%
\global\long\def\stell{\mathrm{stell}}%
\global\long\def\comp{\mathrm{c}}%
\global\long\def\verts{V}%
\global\long\def\edges{E}%
\global\long\def\qq{\ve q}%
\global\long\def\ssv{\ve s}%

\usepackage{natbib}
 \NatBibNumeric
 \def\bibfont{\small}%
 \bibpunct[, ]{[}{]}{,}{n}{}{,}%

\usepackage[colorlinks=true,breaklinks=true,bookmarks=true,urlcolor=blue,
     citecolor=blue,linkcolor=blue,bookmarksopen=false,draft=false]{hyperref}
     
\usepackage{footnotebackref}
     
\def\EMAIL#1{\href{mailto:#1}{#1}}

\TheoremsNumberedThrough     

\EquationsNumberedThrough    

\begin{document}


\RUNAUTHOR{Celaya et al.}

\RUNTITLE{Proximity and flatness bounds for ILPs}
\TITLE{Proximity and flatness bounds for linear integer optimization}

\ARTICLEAUTHORS{%
\AUTHOR{Marcel Celaya}
\AFF{Department of Mathematics, Institute for Operations Research, ETH
Z\"{u}rich, Switzerland, \EMAIL{marcel.celaya@ifor.math.ethz.ch}}
\AUTHOR{Stefan Kuhlmann}
\AFF{Institut f\"{u}r Mathematik, Technische Universit\"{a}t Berlin,
Germany, \EMAIL{kuhlmann@math.tu-berlin.de}}
\AUTHOR{Joseph Paat}
\AFF{Sauder School of Business, University of British Columbia, BC Canada, \EMAIL{joseph.paat@sauder.ubc.ca}}
\AUTHOR{Robert Weismantel}
\AFF{Department of Mathematics, Institute for Operations Research, ETH
Z\"{u}rich, Switzerland, \EMAIL{robert.weismantel@ifor.math.ethz.ch}}
} 

\ABSTRACT{This paper deals with linear integer optimization. 
We develop a technique that can be applied to provide improved upper bounds for two important questions in linear integer optimization.

\smallskip
\begin{itemize}
	\item Proximity bounds: Given an optimal vertex solution for the linear relaxation, how far away is the nearest optimal integer solution (if one exists)?
	\item Flatness bounds: If a polyhedron contains no integer point, what is the smallest number of integer parallel hyperplanes defined by an integral, non-zero, normal vector that intersect the polyhedron?
\end{itemize}

\smallskip
This paper presents a link between these two questions by refining a proof technique that has been recently introduced by the authors. 
A key technical lemma underlying our technique concerns the areas of certain convex polygons in the plane: if a polygon $K\subseteq\R^2$ satisfies $\rot K \subseteq K^{\circ}$, where $\rot$ denotes $90^{\circ}$ counterclockwise rotation and $K^{\circ}$ denotes the polar of $K$, then the area of $K^{\circ}$ is at least 3.
}

%

\maketitle

\section{Introduction.}

Suppose $\AA$ is an integral full-column-rank $m\times n$ matrix. By 
\[
\Delta_{k}\left(\AA\right):=\max\left\{ \left|\det\MM\right|:\MM\text{ is a \ensuremath{k\times k} submatrix of \ensuremath{\AA}}\right\}
\]
we denote the largest absolute $k\times k$ minor of $\AA$.
The polyhedron corresponding to a right hand side $\bb\in\Q^{m}$ is 
\begin{align*}
\poly\left(\AA,\bb\right) & :=\left\{ \xx\in\R^{n}:\ \AA\xx\leq\bb\right\} .\\
\intertext{\text{The linear program corresponding to \ensuremath{\poly\left(\AA,\bb\right)} and an objective vector \ensuremath{\cc\in\Q^{n}} is}}\LP(\AA,\bb,\cc) & :=\max\left\{ \cc^{\top}\xx:\ \xx\in\poly\left(\AA,\bb\right)\right\} ,\\
\intertext{\text{and the corresponding integer linear program is}}\IP(\AA,\bb,\cc) & :=\max\left\{ \cc^{\top}\xx:\ \xx\in\poly\left(\AA,\bb\right)\cap\Z^{n}\right\} .
\end{align*}
Our point of departure is the following foundational result due to Cook, Gerards, Schrijver, and Tardos that has several applications in integer optimization; see \cite{EW2018,FGL2005,JR2018}.
\begin{theorem}[Theorem 1 in~\cite{CGST1986}]\label{thmCooketAl} 
	Let $\bb\in\Q^{m}$ and $\cc\in\Q^{n}$. 
	Let $\xx^{*}$ be an optimal vertex of $\LP(\AA,\bb,\cc)$. 
	If $\IP(\AA,\bb,\cc)$ is feasible, then there exists an optimal solution $\zz^{*}$ such
	that\footnote{Their upper bound is stated as $n\cdot\max\big\{\Delta_{k}(\AA):\ k=1,\ldots,n\big\}$, but their argument actually yields an upper bound of  $n\cdot\Delta_{n-1}\left(\AA\right)$.
		Furthermore, their result holds for any (not necessarily vertex) optimal LP solution $\xx^*$.}
	\[
	\|\xx^{*}-\zz^{*}\|_{\infty}\le  n\cdot\Delta_{n-1}\left(\AA\right).
	\]
\end{theorem}

The technique to prove Theorem \ref{thmCooketAl} has been used to establish proximity bounds involving other data parameters~\cite{L2019} and different norms~\cite{LPSX2020,LPSX2021}.
Furthermore, their result has been extended to derive proximity results for convex separable programs~\cite{G1990,H1990,W1991} (where the bound in Theorem~\ref{thmCooketAl} remains valid), for mixed integer programs~\cite{PWW2020}, and for random integer programs~\cite{OPW2020}.

Lov\'asz~\cite[Section 17.2]{AS1986} and Del Pia and Ma~\cite[Section 4]{DM2021}
identified tuples $(\AA,\bb,\cc)$ such that proximity is arbitrarily
close to the upper bound in Theorem~\ref{thmCooketAl}. However,
their examples crucially rely on the fact that $\bb$ can take arbitrary
rational values. In fact, Lov\'asz's example uses a totally unimodular
matrix $\AA$ while Del Pia and Ma use a unimodular matrix. 
Therefore, if the right hand sides $\bb$ in their examples were to be replaced by the integral rounded down vector $\lfloor\bb\rfloor$, then the polyhedron $\poly(\AA, \lfloor\bb\rfloor)$ would only have integral vertices. From an integer programming perspective, replacing $\bb$ with $\lfloor\bb\rfloor$
is natural as it strengthens the linear relaxation without cutting
off any feasible integer solutions. 

It remains an open question whether Cook et al.'s bound is tight when $\bb\in\Z^{m}$.
Under this assumption, Paat et al.~\cite{PWW2020} conjecture that the true bound is independent of $n$.
This conjecture is supported by various results:
Aliev et al.~\cite{AHO2019} prove that proximity is upper bounded by the largest entry of $\AA$ for knapsack polytopes, Veselov and Chirkov's result~\cite{VC2009} implies a proximity bound of $2$ when $\Delta_{n}(\AA)\le2$, and Aliev et al.~\cite{ACHW2021} prove a bound of $\Delta_{n}(\AA)$ for corner polyhedra.

One of our main results is an improvement on Theorem~\ref{thmCooketAl}.

\begin{theorem}
\label{thm:main_thm_n_2}Let $n\geq2$, $\bb\in\Z^{m}$, and $\cc\in\Q^{n}$.
Let $\xx^{*}$ be an optimal vertex of $\LP(\AA,\bb,\cc)$. If $\IP(\AA,\bb,\cc)$
is feasible, then there exists an optimal solution $\zz^{*}$ such
that
\[
\|\xx^{*}-\zz^{*}\|_{\infty}<\frac{4n+2}{9}\cdot\Delta_{n-1}(\AA).
\]
\end{theorem}

A second equally fundamental question in discrete mathematics is concerned with bounds on flatness of $\poly(\AA,\bb)$ if $\poly(\AA,\bb)$ is lattice-free, i.e., $\poly(\AA,\bb)\cap\Z^n = \emptyset$.
The width of $\poly\left(\AA,\bb\right)$ in direction $\aa \in \R^n\backslash\lbrace 0\rbrace$ is defined by 
\[
w^{\aa}\left(\poly\left(\AA,\bb\right)\right) := \max_{\xx\in \poly\left(\AA,\bb\right)} \aa^\top\xx - \min_{\yy \in \poly\left(\AA,\bb\right)}\aa^\top\yy.
\]
The lattice width is defined by
\[
w\left(\poly\left(\AA,\bb\right)\right) := \min_{\aa\in\Z^n\backslash\lbrace\zero\rbrace} w^{\aa}\left(\poly\left(\AA,\bb\right)\right).
\]
A prominent result regarding the lattice width is due to Khinchine.
\begin{theorem}[\cite{K1948}]\label{thm:khinchine_finiteness}
	Let $\poly\left(\AA,\bb\right)$ be a lattice-free polyhedron. There exists a non-zero vector $\aa\in\Z^n$ such that $w^{\aa}\left(\poly\left(\AA,\bb\right)\right)$ is bounded above by some function depending only on the dimension $n$.
\end{theorem}
The current best upper bound is $\mathcal{O}^*(n^{\frac{4}{3}})$, where $\mathcal{O}^*$ denotes that a polynomial in $\log n$ is omitted; see \cite{R2000}. It is conjectured that the lattice width can be bounded by a function which only depends linearly on $n$.

A variety of algorithms related to integer programs rely on upper bounds on the lattice width of lattice-free polytopes. One famous example is Lenstra's approach to solve the feasibility question of integer linear programs \cite{L1983}.
In order to improve the understanding of the running time of these algorithms with respect to their input, it is a natural task to analyze the lattice width in dependence of other input parameters than $n$.

Gribanov and Veselov presented the first bound on the lattice width of lattice-free polytopes which depends linearly on $n$ and on the least common multiple of all $n\times n$ minors; see \cite{GV2016}. The least common multiple is in the worst case exponentially large in $\Delta_n(\AA)$.
We present a bound that depends linearly on $n$ and linearly on $\Delta_n(\AA)$.
\begin{theorem}\label{thm:facet_width}
	Let $n\geq2$ and $\bb\in\Z^m$ such that $\poly\left(\AA,\bb\right)$ is a full-dimensional lattice-free polyhedron and each row of $\AA$ is facet-defining. Then, there exists a row $\aa$ of $\AA$ such that
	\begin{align*}
		w^\mathcal{\aa}\left(\poly\left(\AA,\bb\right)\right)< \frac{4n+2}{9}\cdot\Delta_{n}(\AA) - 1.
	\end{align*}
\end{theorem}

It is open whether the lattice width of lattice-free polytopes can be bounded solely by $\Delta_n(\AA)$. Some interesting classes of polytopes, where this is the case, are simplices and special pyramids; see \cite{HKW2022}. In \cite{BJ2022}, the authors utilize bounds on the facet width of certain lattice-free polytopes with respect to their minors to construct an algorithm which efficiently enumerates special integer vectors in those polytopes.

On the first glance Theorem~\ref{thm:main_thm_n_2} and \ref{thm:facet_width} have nothing in common. However, we will show that both results follow from a more general result that allows us to establish a bound on the gap of the value between a linear optimization problem and its integer analogue. This applies to arbitrary integral valued objective function vectors. In order to state this result formally, let us introduce the following definition.
\begin{definition}\label{def:delta_alpha}
	Let $\AA\in\Z^{m\times n}$ be full-column-rank matrix. For $\aalpha\in\Z^n$, let
	\begin{align*}
		\Delta^{\aalpha}\left(\AA\right) := \max \left\lbrace \left|\det \left(\begin{array}{l}
			\aalpha^\top \\
			\BB
		\end{array}
		\right)\right| : \BB \text{ is a }(n-1)\times n \text{ submatrix of }\AA\right\rbrace.
	\end{align*}
\end{definition}

\begin{theorem}\label{thm:main_thm_linear_functional}
	Let $\aalpha\in\Z^n\backslash\lbrace \zero\rbrace$, $n\geq2$, $\bb\in\Z^{m}$, and $\cc\in\Q^{n}$.
	Let $\xx^{*}$ be an optimal vertex of $\LP(\AA,\bb,\cc)$. If $\IP(\AA,\bb,\cc)$
	is feasible, then there exists an optimal solution $\zz^{*}$ such
	that
	\[
	\left|\aalpha^\top(\xx^{*}-\zz^{*})\right|<\frac{4n+2}{9}\cdot\Delta^{\aalpha}\left(\AA\right).
	\] 
\end{theorem}

Our proof of Theorem~\ref{thm:main_thm_linear_functional} consists of three major parts. 
First, we apply a dimension reduction technique so that general instances can be reduced to full-dimensional instances in lower dimensional space. This is discussed in Section \ref{secDimensionReduction}. Next, we establish a relationship between Theorem \ref{thm:main_thm_linear_functional} and the volume of a particular polytope associated with the matrix $\AA$ and the vector $\aalpha$. This applies to arbitrary dimensions. In order to give an estimate on the volume, we restrict our attention to low dimensional cases. 
We establish such lower bounds when $n=1$ and $n=2$; see the end of Section~\ref{secDimensionReduction}. Third, we show how these bounds in lower dimensions can be lifted to bounds in higher dimensions; see Section~\ref{secLifting}. Section \ref{sec3dimpolyhedra} is devoted to carrying out the calculations when $n=3$. 

When $n=3$, this particular polytope transforms linearly into a polygon $\qoly^{\circ}\subseteq\R^2$ satisfying $\rot\qoly\subseteq\qoly^{\circ}$, where $\rot$ denotes the $90^\circ$ counterclockwise rotation. In Appendix~\ref{appAreaPolygonPolar}, we show that the area of any such polygon is at least 3. For this we show that it is sufficient to consider the extremal case where $\rot\qoly=\qoly^{\circ}$. Polytopes of this type have been analysed by Jensen in \cite{J2021} (see also \cite{F2020} for the planar case), where they are called \emph{self-polar polytopes}. Next, by repeatedly applying Jensen's \emph{add-and-cut} modification described in \cite{J2021}, we show all minimal-area polygons $\qoly$  satisfying $\rot\qoly=\qoly^{\circ}$ have area 3. We remark that inequalities relating the volume of a polytope with the volume of its polar have long been investigated; their product is the subject of Mahler's conjecture~\cite{M1939B} (see also~\cite[Page 177]{G2007}), and their sum has been studied in the planar case~\cite{F1996}.

Our techniques can be adapted to analyze the special case when the $n\times n$ minors of $\AA$ are contained in $\lbrace 0,\pm k ,\pm 2 k\rbrace$ for some integer $k\geq 1$. A special instance of such a matrix $\AA$ is the case when $\AA$ is \emph{strictly $\Delta_n(\AA)$-modular}, that is, $\AA=\TT\BB$ for a totally unimodular matrix $\TT$ and a square integer matrix $\BB$ with determinant $\Delta_n(\AA)$. In this case, the bounds on proximity and flatness are independent of the dimension, generalizing results of N\"agele, Santiago, and Zenklusen \cite[Theorem 1.4 and 1.5]{NSZ2021}.
\begin{theorem}\label{thm:dim_free_bounds}
	Let the $n\times n$ minors of $\AA$ be contained in $\lbrace 0, \pm k, \pm 2 k\rbrace$ for some integer $k\geq 1$. 
	The following hold:
	\begin{enumerate}
		\item The bound in Theorem \ref{thm:main_thm_n_2} can be sharpened to
		\begin{align*}
			\|\xx^{*}-\zz^{*}\|_{\infty}\leq\max\lbrace\Delta_{n-1}(\AA),\Delta_n(\AA)\rbrace - 1
		\end{align*}
		and 
		\begin{align*}
			\|\xx^{*}-\zz^{*}\|_{\infty}<\Delta_{n-1}(\AA).
		\end{align*}
		\item The bound in Theorem \ref{thm:facet_width} can be sharpened to 
		\begin{align*}
			w^\mathcal{\aa}\left(\poly\left(\AA,\bb\right)\right)\leq\Delta_{n}(\AA) - 2.
		\end{align*}
	\end{enumerate}
\end{theorem}
The proof of this theorem is given in Section \ref{secDimensionFree}.


\medskip

\begin{remark}
This manuscript builds upon the work carried out by the authors in \cite{CKPW2022}. 
All of the results in this paper are strict improvements of the results in~\cite{CKPW2022}, with the main new contributions being the improved constant in Theorem~\ref{thm:main_thm_n_2}, the flatness result of Theorem~\ref{thm:facet_width}, and the extension of Theorem~\ref{thm:dim_free_bounds} to the $\{0,\pm k,\pm 2k\}$-setting. \hfill$\diamond$
\end{remark}

\section{Basic Definitions and Notation.}\label{ssecPrelim_not}
Here we outline the key objects and parameters used in the paper.

Let $\AA\in\Z^{m\times n}$ be a full-column-rank matrix, and $\bb\in\Z^{m}$ be such that $\poly(\AA, \bb) \cap \Z^n \neq \emptyset$.
For $I \subseteq [m] := \{1, \ldots, m\}$, we use $\AA_I$ and $\bb_I$ to denote the rows of $\AA$ and $\bb$ indexed by $I$.
If $I = \{i\}$, then we write $\aa_i^\top  := \AA_I$.
We use $\zero$ and $\one$ to denote the all zero and all one vector (in appropriate dimension).
For a polyhedron $\mathcal{Q} \subseteq \R^n$, the dimension of $\mathcal{Q}$ is the dimension of the linear span of $\mathcal{Q}$ and is denoted by $\dim \mathcal{Q}$. We also define, for $I\subseteq[m]$,
\[
\gcd\AA_{I} := \gcd \left\{\left|\det \MM \right|:\ \MM~\text{is a}~\text{rank}(\AA_I) \times \text{rank}(\AA_I)~\text{submatrix of}~\AA_I\right\},
\]
with $\gcd\AA_{\emptyset}=1$. In the case when $\poly\cap\Z^n=\{\zero\}$, bounding proximity is equivalent to bounding 
\begin{equation}\label{eqProxWidth}
\max_{\xx \in \poly(\AA, \bb)}\ \|\xx\|_{\infty} = \max_{\aalpha \in \{\pm\ee_1, \ldots, \pm \ee_n\} }\ \max \left\{\aalpha^\top \xx:\  \xx \in \poly(\AA, \bb)\right\},
\end{equation}
where $\ee_1, \ldots, \ee_n \in \Z^n$ are the standard unit vectors. As we shall see in the proof of Theorem~\ref{thm:main_thm_linear_functional}, the general case then follows from this case.
In light of this, we analyze the maximum of an arbitrary linear form $\aalpha^\top \xx$ over $\poly(\AA, \bb)$ for $\aalpha \in \Z^n$.

We provide non-trivial bounds on the maximum of these linear forms for small values of $n$; see Section~\ref{secDimensionReduction} and \ref{sec3dimpolyhedra}.
In order to lift low dimensional results to higher dimensions (see Section~\ref{secLifting}), we consider slices of $\poly\left(\AA,\bb\right)$ through the origin induced by rows of $\AA$. 
Given $I\subseteq\left[m\right]$ such that $|I| \le n-1$ and $\text{rank}~ \AA_I = |I|$, define
\[
\poly_{I}\left(\AA,\bb\right):=\poly\left(\AA,\bb\right)\cap\ker\AA_{I}.
\]
We specify $\ker\AA_{\emptyset}=\R^n$, so that $\poly_{\emptyset}(\AA, \bb) = \poly(\AA, \bb)$.
The bounds that we provide on $\aalpha^\top \xx$ are given in terms of the parameter
\[
\Delta^{\aalpha}_{I}\left(\AA\right):=\frac{1}{\gcd\AA_{I}}\cdot\max\left\{ \left|\det\left(\begin{array}{c}
\objfn^\top\\
\AA_{K}
\end{array}\right)\right|:I\subseteq K\subseteq\left[m\right],\;\left|K\right|=n-1\right\}.
\]
Observe that $\Delta^{\aalpha}\left(\AA\right)=\Delta^{\aalpha}_{\emptyset}(\AA)$. 
In particular, we define $\prox_{I}\left(\AA,\bb,\objfn\right)$ to be the number satisfying
\begin{equation}\label{eqPseudoProx}
\max_{\xx\in\poly_{I}\left(\AA,\bb\right)}\objfn^\top\xx=\prox_{I}\left(\AA,\bb,\objfn\right)\Delta^{\aalpha}_{I}\left(\AA\right).
\end{equation}
Maximizing over all $I\subseteq\left[m\right]$ such that $\poly_{I}\left(\AA,\bb\right)$ has a fixed dimension $d$, define
\[
\prox_{d}\left(\AA,\bb,\objfn\right):=\max_{I:\dim\poly_{I}\left(\AA,\bb\right)=d}\prox_{I}\left(\AA,\bb,\objfn\right).
\]
Equation~\eqref{eqPseudoProx} looks similar to the bound we seek. 
However, $\Delta^{\aalpha}_{I}\left(\AA\right)$ depends on $\objfn$, whereas our main result (Theorem~\ref{thm:main_thm_n_2}) only depends on $\Delta_{n-1}(\AA)$.
Later (see Section~\ref{secProofOfMainTheorem}), we will substitute $\pm \ee_1, \ldots, \pm\ee_n$ in for $\aalpha$ as in~\eqref{eqProxWidth}.
We also want to consider $I = \emptyset$ because $\poly_{\emptyset}(\AA, \bb)  = \poly(\AA, \bb)$ by definition. Note that when $\aalpha$ is a unit vector, then $\Delta^{\aalpha}_{I}\left(\AA\right)$ is a lower bound for $\Delta_{n-1}(\AA)$.
Another important object for us is the following cone.
For $\xx^{*}\in\R^{n}$, define
\begin{align*}
\ccone\left(\AA,\xx^{*}\right):=\left\{\xx\in\R^{n}:
\begin{array}{rl}
\sign\left(\aa_{i}^\top\xx^{*}\right)\cdot\aa_{i}^\top\xx & \geq0~\forall~ i \in [m]~\text{such that}~\aa_{i}^\top\xx^{*}\neq0\\[.1 cm]
\aa_{i}^\top\xx & =0~\forall~ i \in [m]~\text{such that}~\aa_{i}^\top\xx^{*}=0
\end{array}\right\}.
\end{align*}
The cone $\ccone\left(\AA,\xx^{*}\right)$ serves as a key ingredient in the proof of Theorem~\ref{thmCooketAl} in \cite{CGST1986}. 
We also define the polytope
\[
\spindle\left(\AA,\xx^{*}\right):=\ccone\left(\AA,\xx^{*}\right)\cap\left(\xx^{*}-\ccone\left(\AA,\xx^{*}\right)\right).
\]
One checks that if $\xx^{*}\in\poly\left(\AA,\bb\right)$, then $\spindle\left(\AA,\xx^{*}\right)\subseteq\poly\left(\AA,\bb\right)$.
Moreover, if $\yy^{*}\in\spindle\left(\AA,\xx^{*}\right)$ then $\spindle\left(\AA,\yy^{*}\right)\subseteq\spindle\left(\AA,\xx^{*}\right)$.
Polytopes of this form, namely, ones in which every facet is incident to one of two distinguished vertices, known as \emph{spindles}, were used in \cite{S2012} to construct counterexamples to the Hirsch conjecture.

We often fix $\AA\in\Z^{m\times n}$ and $\bb\in\Z^{m}$. 
Thus, if the dependence on $\AA$ and $\bb$ is clear from the context, we abbreviate $\poly_{I}$ for $\poly_I(\AA,\bb)$, $\Delta^{\aalpha}_{I}$ for $\Delta^{\aalpha}_{I}\left(\AA\right)$, $\prox(\aalpha)$ for $\prox(\AA,\bb,\aalpha)$, $\spindle(\xx^*)$ for $\spindle(\AA,\xx^*)$ and so on.

\section{Dimension Reduction and Further Preliminaries.}\label{secDimensionReduction}

A useful fact for us is that we only need to consider the case when $\dim \poly = n$, by replacing a not-necessarily full-dimensional instance with an equivalent full-dimensional instance in a lower-dimensional space. This construction is outlined below.
%

\begin{lemma}
\label{lem:lifting_lemma}
Let $\aalpha\in\Z^n$ such that $\max\{\aalpha^\top\xx:\xx\in\poly\}$ is attained and is finite. Assume $I\subseteq\left[m\right]$ determines a linearly independent subset of the rows of $\AA$ such that the linear span of $\poly_{I}$ is $\ker\AA_{I}$, which has dimension $d$. 
Then there exists a linear isomorphism $\ker\AA_I\rightarrow \R^d$ given by $\xx\mapsto\PP\xx$ where $\PP\in\Z^{d\times n}$, which maps $\ker\AA_I\cap\Z^n$ onto $\Z^d$ and maps $\poly_I(\AA,\bb)$ onto $\poly(\AAhat,\hat{\bb})$ for some $\AAhat\in\Z^{\left(m-n+d\right)\times d}$, $\hat{\bb}\in\Z^{m-n+d}$, and satisfies
\[
\prox_{I}\left(\AA,\bb,\objfn\right)=\prox_d\bigl(\AAhat,\hat{\bb},\hat{\objfn}\bigr)
\]
where $\hat{\objfn}\in\Z^d$ is the unique vector satisfying $\hat{\objfn}^\top\PP=\objfn^\top$.

\end{lemma}

\proof{Proof.}
Without loss of generality, suppose $I = [n-d]$. 
Set $J := [n-d]$, $\bar{J} := \{n-d+1, \ldots, n\}$, and $\bar{I} := \{n-d+1, \ldots, m\}$.
Choose a unimodular matrix $\UU\in\Z^{n\times n}$ (e.g., via the Hermite Normal Form of $\AA_{[n]}$) such that 
\[
\AA\UU=\begin{pmatrix}\left(\AA\UU\right)_{I,J} & \zero\\
\left(\AA\UU\right)_{\bar{I},J} & \left(\AA\UU\right)_{\bar{I},\bar{J}}
\end{pmatrix}
\]
with $\left(\AA\UU\right)_{I,J}$ square and invertible. 

Set 
\(
\AAhat :=\left(\AA\UU\right)_{\bar{I},\bar{J}},
\)
\(
\hat{\bb} := \bb_{\bar{I}},
\)
and
\(
\hat{\objfn}^\top :=\left(\objfn^\top\UU\right)_{\bar{J}}.
\)
For $\xx \in \ker \AA_I$, we have 
\[
\zero = \AA_I \xx = \AA_I\UU \UU^{-1} \xx = [(\AA \UU)_{I,J} ~ \zero]\ \UU^{-1} \xx = (\AA\UU)_{I,J}(\UU^{-1}\xx)_J.
\]
Thus, $(\UU^{-1}\xx)_J = \zero$.
Hence, the map $\xx\mapsto\left(\UU^{-1}\xx\right)_{\bar{J}}$ is a linear isomorphism from $\ker\AA_{I}$ to $\R^{|\bar{J}|}=\R^{d}$, which restricts to a lattice isomorphism from $\ker\AA_{I}\cap\Z^{n}$ to $\Z^{d}$ and maps $\poly_{I}(\AA, \bb)$ to $\poly\bigl(\AAhat,\hat{\bb}\bigr)$.
It follows that $\poly\bigl(\AAhat,\hat{\bb}\bigr)\cap\Z^{d}=\left\{ \zero\right\} $.
For $\xx\in\ker\AA_{I}$, the equation $\left(\UU^{-1}\xx\right)_{J}=\zero$ implies that
\begin{equation}
\objfn^\top\xx=\objfn^\top\UU\UU^{-1}\xx=\hat{\objfn}^\top\left(\UU^{-1}\xx\right)_{\bar{J}}.\label{eq:numerator}
\end{equation}
Moreover, if $K\subseteq\bar{I}$ with $\left|K\right|=d-1$, then
\begin{align*}
\left|\det\left(\begin{array}{c}
\hat{\objfn}^\top\\
\AAhat_{K}
\end{array}\right)\right| 
& =\left|\det\left(\begin{array}{r@{\hskip 0 cm}r@{\hskip 0 cm}l}
\left(\right.&\objfn^\top\UU &\left.\right)_{\bar{J}}\\
\left(\right.&\AA\UU&\left.\right)_{K,\bar{J}}
\end{array}\right)\right|\\
 & =\frac{1}{\bigl|\det\left(\AA\UU\right)_{I,J}\bigr|}\cdot
 \left|\det\left(
 \begin{array}{r@{\hskip 0 cm}r@{\hskip 0 cm}lr@{\hskip 0 cm}r@{\hskip 0 cm}l}
 \left(\right.&\AA\UU&\left.\right)_{I,J} & &\zero~~&\\
\left(\right.&\objfn^\top\UU&\left.\right)_{J} & \left(\right.&\objfn^\top\UU&\left.\right)_{\bar{J}}\\
\left(\right.&\AA\UU&\left.\right)_{K,J} & \left(\right.&\AA\UU&\left.\right)_{K,\bar{J}}
\end{array}
\right)\right|\\
 & =\frac{1}{\gcd\AA_{I}}\cdot\left|\det\left(\begin{array}{c}
\objfn^\top\\
\AA_{I\cup K}
\end{array}\right)\right|,
\end{align*}
where we have used 
\(
\bigl|\det\left(\AA\UU\right)_{I,J}\bigr| = \gcd (\AA\UU)_I = \gcd \AA_I \UU = \gcd \AA_I.
\)
Taking the maximum over all such $K$, we get
\begin{equation}
\Delta^{\hat{\objfn}}(\AAhat)=\Delta^{\aalpha}_{I}\left(\AA\right).\label{eq:denominator}
\end{equation}
Putting~\eqref{eq:numerator} and~\eqref{eq:denominator} together,
we get
\[
\prox_d\bigl(\AAhat,\hat{\bb},\hat{\objfn}\bigr)
=\max_{\yy\in\poly\left(\AAhat,\hat{\bb}\right)}\ \frac{\hat{\objfn}^\top\yy}{\Delta^{\hat{\objfn}}(\AAhat)}
=\max_{\xx\in\poly_{I}(\AA, \bb)}\ \frac{\objfn^\top\xx}{\Delta^{\aalpha}_{I}\left(\AA\right)}
=\prox_{I}\left(\AA,\bb,\objfn\right).\tag*{\halmos}
\]
\endproof
Next, we present a general relationship between the volume of polyhedra associated with the matrix $\AA$, $\aalpha$, and $\prox_{n}(\aalpha)$.

Define the polyhedron
\[
\poly_{\objfn}:=\left\{ \xx\in\R^{n}:\left|\AA\xx\right|\leq\one,\;\objfn^\top\xx=0\right\} .
\]
This is an $\left(n-1\right)$-dimensional polyhedron, which is bounded since $\AA$ has full-column-rank by assumption.
We use $\vol_i(\cdot)$ to denote the $i$-dimensional Lebesgue measure.

\begin{lemma}
	\label{lem:prox_volume}
	Let $\aalpha\in\Z^n$ be non-zero. Assume $ \dim \poly = n$ and $\poly\cap\Z^n=\{\zero\}$.
	%
	Then
	\[
	\prox_{n}(\aalpha)<\frac{2^{n-1}\norm{\objfn}_{2}}{\vol_{n-1}\left(\poly_{\objfn}\right)\Delta^{\aalpha}}.
	\]
\end{lemma}

\proof{Proof.}
	Recall $\poly = \poly(\AA, \bb)$.
	Let $\xx^{*}\in\poly$  attain the maximum of 
	\[
	\prox_{n}(\aalpha)=\max_{\xx\in\poly}\frac{\objfn^\top\xx}{\Delta^{\aalpha}},
	\]
	which we assume is positive without loss of generality. 
	Define the polytope
	\[
	\qoly\left(\xx^{*}\right):=\poly_{\objfn}+\left[-\xx^{*},\xx^{*}\right],
	\]
	which is $\zero$-symmetric and full-dimensional in $\R^{n}$. 
	Observe that
	\[
	\vol_{n}\left(\qoly\left(\xx^{*}\right)\right)=\frac{2\prox_{n}(\aalpha)\Delta^{\aalpha}}{\norm{\objfn}_{2}}\cdot\vol_{n-1}\left(\poly_{\objfn}\right).
	\]
	All integer points not in $\qoly(\xx^*)$ are a positive distance away from $\qoly(\xx^*)$, hence there exists $\delta>0$ such that $\qoly((1+\delta)\xx^*)$ and $\qoly(\xx^*)$ contain precisely the same set of integer points. This choice of $\delta$ uniquely determines $\varepsilon>0$ for which 
	\[
	\qoly'\left(\xx^{*}\right):=\left(1-\varepsilon\right)\qoly\left(\left(1+\delta\right)\xx^{*}\right)
	\]
	has the same $n$-dimensional volume as $\qoly(\xx^*)$, and furthermore
	\[
	\qoly'\left(\xx^{*}\right)\cap\Z^{n}\subseteq\qoly\left(\xx^{*}\right)\cap\Z^{n}.
	\]
	Assume to the contrary that $\vol_{n}\left(\qoly\left(\xx^{*}\right)\right)\geq2^{n}$.
	By Minkowski's convex body theorem, there exists $\zz^{*} \in \qoly(\xx^*)\cap\qoly'(\xx^*)\cap\Z^n\setminus\{\zero\}$ by the above inclusion. Therefore, with respect to the vector space decomposition of $\R^n$ into the line $\R\cdot\xx^*$ and the hyperplane $\aalpha^\top\xx=0$, the vector $\zz^*$ decomposes uniquely as $\zz^{*}=\lambda\xx^{*}+\left(\zz^{*}-\lambda\xx^{*}\right)$ with $\lambda\in[0,1]$ and $\zz^{*}-\lambda\xx^{*}\in\left( 1 - \varepsilon \right)\poly_{\aalpha}$.
	Hence,
	\[
	\left|\AA\left(\zz^{*}-\lambda\xx^{*}\right)\right|\leq\left(1-\varepsilon\right)\one.
	\]
	As $\poly\cap\Z^{n}=\left\{ \zero\right\} $ and $\zz^{*} \neq \zero$, there exists some row $\aa_{j}^\top$ of $\AA$ such that $\aa_{j}^\top\zz^{*}\geq\bb_{j}+1$. 
	Since $\xx^{*}\in\poly(\AA, \bb)$, we also have $\aa_{j}^\top\xx^{*}\leq\bb_{j}$. 
	Thus, we get
	\begin{align*}
		\bb_{j}+1\leq\aa_{j}^\top\zz^{*} & =\aa_{j}^\top\left(\lambda\xx^{*}\right)+\aa_{j}^\top\left(\zz^{*}-\lambda\xx^{*}\right)\leq\lambda\bb_{j}+\left(1-\varepsilon\right)<\bb_{j}+1.
	\end{align*}
	This is a contradiction.
	Hence,
	\[
	\frac{2\prox_{n}(\aalpha)\Delta^{\aalpha}}{\norm{\objfn}_{2}}\cdot\vol_{n-1}\left(\poly_{\objfn}\right)=\vol_{n}\left(\qoly\left(\xx^{*}\right)\right)<2^{n}.
	\]
	Rearranging yields the desired inequality.\hfill\halmos
\endproof

\medskip

\begin{remark}
	Integrality of $\bb$, which is the key assumption of this paper, is used above in the assertion $\aa_{j}^\top\zz^{*}\geq\bb_{j}+1$.
	If $\bb$ were not integral, then we would only be able to assert that $\aa_{j}^\top\zz^{*}\geq\lceil\bb_{j}\rceil$, which is not sufficient to complete the proof. \hfill $\diamond$
\end{remark}

\medskip

A final step in this Section is to establish basic bounds on $\prox_{1}(\aalpha)$ and $\prox_{2}(\aalpha)$. 

%
%

\begin{lemma}
\label{lem:prox_dim_1_2_3}Let $\aalpha\in\Z^n$ be non-zero. Suppose $\poly\cap\Z^n=\{\zero\}$. Then $\prox_{1}(\aalpha)<1$ and  $\prox_{2}(\aalpha)<1$.
\end{lemma}
\proof{Proof.}
By Lemma \ref{lem:lifting_lemma} we may assume $\poly$ is full-dimensional. 
If $n=1$, then $\poly\left(\AA,\bb\right)$ is contained in the open interval $\left(-1,1\right)$, which immediately implies $\prox_{1}(\aalpha)<1$.
If $n=2$, then the polytope $\poly_{\objfn}$ is an origin-symmetric line segment $\left[-\yy^{*},\yy^{*}\right]$, where $\yy^{*}\in\R^{2}$ satisfies $\objfn^\top\yy^{*}  =0$ and $\aa_{j}^\top\yy^{*}  =1$ for some $j\in\left[m\right]$. 
Hence
\[
\vol_{1}\left(\poly_{\objfn}\right)=2\norm{\yy^{*}}_{2}=\frac{2\norm{\objfn}_{2}}{\left|\det\left(\objfn,\aa_{j}\right)\right|}.
\]
Applying Lemma \ref{lem:prox_volume}, we get
\[
\prox_{2}(\aalpha)<\frac{2\norm{\objfn}_{2}}{\vol_{1}\left(\poly_{\objfn}\right)\Delta^{\aalpha}}=\frac{\left|\det\left(\objfn,\aa_{j}\right)\right|}{\Delta^{\aalpha}}\leq1. \tag*{\halmos}
\]
\endproof

\section{An Analysis for 3-Dimensional Polyhedra.}\label{sec3dimpolyhedra}

Recall the definition of the polar $\qoly^{\circ}$ of a non-empty compact convex set $\qoly\subseteq\R^{2}$:
\[
\qoly^{\circ}:=\left\{ \xx\in\R^{2}:\yy^\top\xx\leq1\text{ for all }\yy\in \qoly\right\} .
\]
Also recall that $\rot:\R^2 \to \R^2$ denotes the $90^{\circ}$ counterclockwise rotation in $\R^{2}$. Our bound for $\prox_{3}(\aalpha)$ relies on the following result, which is proved in Appendix~\ref{appAreaPolygonPolar}:

\begin{lemma}
\label{lem:polar_area_lower_bound}
Suppose $\qoly$ is a polygon such that $\rot \qoly \subseteq \qoly^{\circ}$. Then $\vol_2(\qoly^{\circ})\geq 3$.
\end{lemma}

\begin{lemma}
\label{lem:prox_dim_3}Let $\aalpha\in\Z^3$ be non-zero. Suppose $\poly\cap\Z^3=\{\zero\}$. Then $\prox_{3}(\aalpha)<4/3$.
\end{lemma}
\proof{Proof.}
By Lemma \ref{lem:lifting_lemma} we may assume $\poly$ is full-dimensional. 
Choose $I\subseteq\left[m\right]$ with $\left|I\right|=2$ such that
\[
\BB:=\begin{pmatrix}\objfn^\top\\
\AA_{I}
\end{pmatrix}
\]
satisfies $\left|\det\BB\right|=\Delta^{\objfn}$. 
Let $\AA'$ denote the last two columns of $\AA\BB^{-1}$, and enumerate the rows of $\AA'$ as $\aa_{1}',\ldots,\aa_{m}'$. Let $\qoly$ denote the convex hull of these rows and their negatives. Then
\begin{align*}
\BB\cdot\poly_{\objfn} & =\left\{ 0\right\} \times\qoly^{\circ}.
\end{align*}
Since $\poly_{\aalpha}$ is bounded, so is $\qoly^{\circ}$. Observe that
\(
\rot\qoly\subseteq\qoly^{\circ}.
\)
Indeed, for each pair $\{i,j\}\subseteq\left[m\right]$, we have
\[
\bigl|(\rot\aa_{i}')^\top\aa_{j}'\bigr|=\bigl|\det(\aa_{i}',\aa_{j}')\bigr|=\frac{\bigl|\det(\objfn,\aa_{i},\aa_{j})\bigr|}{\bigl|\det\BB\bigr|}\leq1.
\]
Hence, by Lemma~\ref{lem:polar_area_lower_bound}, we get $\vol_{2}\left(\qoly^{\circ}\right)\geq 3$. We have
\[
\vol_{2}\left(\qoly^{\circ}\right)=\frac{\left|\det\BB\right|}{\norm{\objfn}_{2}}\cdot\vol_{2}\left(\poly_{\objfn}\right),
\]
and so by Lemma~\ref{lem:prox_volume}, we get
\begin{align*}
\prox_{3}(\aalpha)<\frac{4\norm{\objfn}_{2}}{\vol_{2}\left(\poly_{\objfn}\right)\Delta^{\aalpha}}=\frac{4}{\vol_{2}\left(\qoly^{\circ}\right)}\leq\frac{4}{3}.\tag*{\halmos}
\end{align*}
\endproof

\section{Lifting Low Dimensional Results to Higher Dimensions.}\label{secLifting}

The next step is to prove Theorem~\ref{thm:main_thm_n_2} by showing how results for low dimensional polytopes can be used to derive results for higher dimensional polytopes.

\begin{lemma}
\label{lem:spindle_down_to_zero}
Let $\xx^*\in\R^n$, and let $d=\dim\spindle\left(\xx^{*}\right)$. Let $\yy^{*}\in\spindle\left(\xx^{*}\right)$, let $k:=\dim \spindle\left(\yy^{*}\right)$, and fix $d \in \{1, \ldots, k\}$.
There exists a $d$-face of $\spindle\left(\yy^{*}\right)$ incident to $\yy^{*}$ that intersects some $\left(k-d\right)$-face of $\spindle\left(\yy^{*}\right)$ incident to $\zero$.
\end{lemma}

\proof{Proof.}
Let $I\subseteq\left[m\right]$ index the components $i$ such that $\aa_{i}^\top\yy^{*}\neq0$. For $i\in I$ let $\hat{\aa}_{i}=\sign\left(\aa_{i}^\top\yy^{*}\right)\cdot\aa_{i}$.
The spindle $\spindle\left(\yy^{*}\right)$ can be written as
\begin{align*}
\spindle\left(\yy^{*}\right) = \left\{\xx \in \R^n:\ \zero\leq\hat{\aa}_{i}^\top\xx \leq\hat{\aa}_{i}^\top\yy^{*} ~\forall ~i \in I~\text{and}~\aa_{i}^\top\xx  =0 ~\forall ~ i \not \in I\right\}.
\end{align*}
The constraints are indexed by the disjoint union $I_{\zero}\cup I_{\yy^{*}}\cup\bar{I}$, where $I_{\zero}$ and $I_{\yy^{*}}$ denote the two copies of $I$ indexing constraints tight at $\zero$ and at $\yy^{*}$, respectively.
Let $J_{0},J_{1},\ldots,J_{r}$ be a sequence of feasible bases of this system, with corresponding basic feasible solutions $\zero=\yy^{(0)},\yy^{(1)},\ldots,\yy^{(r)}=\yy^{*}$ such that for each $i<r$, the symmetric difference of $J_{i+1}$ and $J_{i}$ is a 2-element subset of $I_{\zero}\cup I_{\yy^{*}}$.
We have $|J_{0}\cap I_{\yy^{*}}|=0$ and $|J_{r}\cap I_{\yy^{*}}|=k$, and $|J_{i+1}\backslash J_{i}|=1 $ for each $i<r$. 
It follows that there must exist some $\ell$ such that $|J_{\ell}\cap I_{\yy^{*}}|=k-d$.
Since we always have $|J_{i}\cap(I_{\zero}\cup I_{\yy^{*}})|=k$ for every choice of $i$, we also get $|J_{\ell}\cap I_{\zero}|=d$. 

The basic feasible solution $\yy^{(\ell)}$ associated to $J_{\ell}$ is a vertex of the face of $\spindle\left(\yy^{*}\right)$ obtained by making the constraints of $J_{\ell}\cap I_{\yy^{*}}$ tight. 
It is also a vertex of the face of $\spindle\left(\yy^{*}\right)$ obtained by making the constraints of $J_{\ell}\cap I_{\zero}$ tight. 
These faces are contained in a $d$-face and a $\left(k-d\right)$-face, respectively.
\hfill\halmos
\endproof

Lemma~\ref{lem:spindle_down_to_zero} will be used to create a path from one vertex of a spindle to another by traveling over $d$ dimensional faces.
In the next result, we apply $d$ dimensional results to each $d$ dimensional face that we travel over.
This generalizes the proof of Cook et al., which can be interpreted as walking along edges of a spindle.

\begin{lemma}
\label{lem:proximity_template}
Let $\aalpha\in\Z^n$ be non-zero. Let $\dim\poly =: d = \sum_{i=0}^k d_i$ where each $d_{i}$ is a positive integer. 
Then
\[
\prox_{d}(\aalpha)\leq\textstyle \sum_{i=0}^k \prox_{d_i}(\aalpha).
\]
\end{lemma}

\proof{Proof.}
In this proof, we suppress in our notation dependence on $\aalpha$. Let $\xx^{*}$ maximize $\aalpha^\top\xx$ over $\poly$.
%
Build a sequence $\xx^{*}=:\xx_{0}^{*},\xx_{1}^{*},\ldots,\xx_{t}^{*}:=\zero$ of points inductively as follows. 
Assume $i\geq0$ and $\xx_{0}^{*},\ldots,\xx_{i}^{*}$ have been determined already.
If both 
\begin{equation}
i\leq k~\text{ and }~d_{i}<\dim\spindle\left(\xx_{i}^{*}\right),\label{eq:two_conditions}
\end{equation}
then we use Lemma~\ref{lem:spindle_down_to_zero} to choose a vertex $\xx_{i+1}^{*}$ of $\spindle\left(\xx_{i}^{*}\right)$ that is incident to both a $d_{i}$-dimensional face $F_{i}$ of $\spindle\left(\xx_{i}^{*}\right)$ containing $\xx_{i}^{*}$, as well as a $\left(\dim\spindle\left(\xx_{i}^{*}\right)-d_{i}\right)$-dimensional face $G_{i}$ of $\spindle\left(\xx_{i}^{*}\right)$ containing $\zero$.
Otherwise, if~\eqref{eq:two_conditions} fails, then we set $F_{i}=\spindle\left(\xx_{i}^{*}\right)$ and $\xx_{i+1}^{*}=\zero$, and we terminate the sequence by setting $t=i+1$.

Let $i \in \{0, \ldots, t-2\}$.
%
We show $\xx_{i+1}^*\neq\zero$. If not, then $F_i$ contains both $\zero$ and $\xx_{i}^*$. But the only face of $\spindle\left(\xx_{i}^{*}\right)$ containing $\zero$ and $\xx_{i}^*$ is $\spindle\left(\xx_{i}^{*}\right)$ itself. One can see this by observing that the centre of symmetry of the centrally symmetric spindle $\spindle\left(\xx_{i}^{*}\right)$ is $\sfrac{1}{2}\cdot\xx_i^*$. But this contradicts the fact that $G_i$ has positive dimension by \eqref{eq:two_conditions}.
Thus, $\xx_{i+1}^{*}$ is non-zero, which implies
\begin{equation}\label{eqBigDim}
\dim\spindle\left(\xx_{i+1}^{*}\right)\geq1.
\end{equation}
Moreover, as both $G_{i}$ and $\spindle\left(\xx_{i+1}^{*}\right)$ are contained in the affine (equivalently, linear) span of $G_{i}$, we must have
\begin{equation}\label{eqStepThrough}
\dim\spindle\left(\xx_{i+1}^{*}\right)\leq\dim G_{i}=\dim\spindle\left(\xx_{i}^{*}\right)-d_{i}.
\end{equation}
Applying~\eqref{eqBigDim} and then~\eqref{eqStepThrough} sequentially with $s \in \{t-2,t-3,\ldots,0\}$, we have
\[
1\leq\dim\spindle\left(\xx_{t-1}^{*}\right) 
  \leq\dim\spindle\left(\xx_{0}^{*}\right)-\sum_{s=0}^{t-2}\ d_{s} 
  \leq d-\sum_{s=0}^{t-2}\ d_{s},
\]
which is to say
\(
d=\sum_{s=0}^{k} d_s > \sum_{s=0}^{t-2}d_s.
\)
It follows that $t-1\leq k$. 

Suppose $I\subseteq\left[m\right]$ indexes linearly independent rows of $\AA$ such that $\prox_d=\prox_I$, so that in particular $\ker\AA_{I}$ is the linear span of $\poly$. Let $i \in \{0, \ldots, t-1\}$.
We have that $\xx_{i}^{*}-F_{i}$ is a face of $\spindle\left(\xx_{i}^{*}\right)$ containing $\zero$.
Choose an index set $I_i$, where $I\subseteq I_{i}\subseteq\left[m\right]$, such that the rows of $\AA_{I_i}$ are linearly independent and $\ker\AA_{I_{i}}$ is the linear span of $\xx_{i}^{*}-F_{i}$. 
We have
\[
\objfn^\top\left(\xx_{i}^{*}-\xx_{i+1}^{*}\right)\leq\max_{\xx\in\xx_{i}^{*}-F_{i}}\objfn^\top\xx\leq\max_{\xx\in\poly_{I_{i}}}\objfn^\top\xx\leq\prox_{I_{i}}\Delta_{I_{i}}.
\]
If $i<t-1$, then since $F_i$ is a $d_i$-dimensional face, we have
\(
\prox_{I_{i}}\Delta_{I_{i}}\leq \prox_{d_{i}}\Delta_{I}
\)
for $i\in \{0,\ldots,t-2\}$.
Otherwise $i=t-1$, in which case one of the inequalities in~\eqref{eq:two_conditions} fails. We have established that $t-1\leq k$, thus
\begin{equation*}
d_{t-1}\geq\dim\spindle\left(\xx_{t-1}^{*}\right)=\dim F_{t-1}.
\end{equation*}
and hence
\(
\prox_{I_{t-1}}\Delta_{I_{t-1}}\leq \prox_{d_{t-1}}\Delta_{I}
\).
Putting these all together we get
\[
\Delta_I\cdot\prox_d = \objfn^\top\xx^* =
\sum_{i=0}^{t-1}\ \objfn^\top\left(\xx_{i}^{*}-\xx_{i+1}^{*}\right)\ \leq\ \sum_{i=0}^{t-1}\ \prox_{I_{i}}\Delta_{I_{i}}\
\leq\ \Delta_I \cdot \sum_{i=0}^k\prox_{d_i}. \tag*{\halmos}
\]
\endproof

\section{Proof of Theorem \ref{thm:main_thm_linear_functional} and \ref{thm:main_thm_n_2}.}\label{secProofOfMainTheorem}
The first step of the proof of \ref{thm:main_thm_linear_functional} is the following reduction which turns out to be useful in later sections.
\begin{lemma}\label{lem:only_one_integer_vector}
	Given $\poly\left(\AA,\bb\right)$ and $\xx^*$, an optimal vertex of $\LP(\AA, \bb, \cc)$, then there exist $\zz^*$, an optimal solution of $\IP(\AA, \bb, \cc)$, an integral matrix $\overline{\AA}$, and an integral vector $\overline{\bb}$ such that $\poly\left(\overline{\AA},\overline{\bb}\right)\subseteq \poly\left(\AA,\bb\right)$, $\xx^*\in\poly\left(\overline{\AA},\overline{\bb}\right)$ is a vertex, and $\poly\left(\overline{\AA},\overline{\bb}\right)\cap\Z^n = \lbrace \zz^*\rbrace$, where the rows of $\overline{\AA}$ consists of rows of $\AA$ and their negatives.
\end{lemma}

\proof{Proof.}
	By LP duality, there exists an optimal LP basis $I^* \subseteq[m]$, i.e., $\xx^* = \AA_{I^*}^{-1} \bb_{I^*}^{\phantom{-1}}$, and a vector $\yy \in \R^{I^*}_{\ge 0}$ that satisfies $\cc^\top=\yy^\top\AA_{I^*}$.
	The polytope $\poly(\overline{\AA}, \overline{\bb}):= \{\xx \in \poly(\AA, \bb):\ \AA_{I^*} \xx \ge \AA_{I^*} \zz^*\}$ contains $\xx^*$ and $\zz^*$ and $\Delta_k(\overline{\AA}) = \Delta_k(\AA)$ for all $k \in [n]$.
	Any integer vector $\ww^* \in \poly(\overline{\AA}, \overline{\bb}) \setminus \{\zz^*\}$ is also an optimal solution to $\IP(\AA, \bb, \cc)$ because $\cc^\top \ww^* = \yy^\top (\AA_{I^*}\ww^*) \ge \yy^\top (\AA_{I^*}\zz^*) = \cc^\top \zz^*$, and $\AA_{I^*} \ww^* \ge \AA_{I^*} \zz^*$ with at least one of the $n$ inequalities satisfied strictly because $I^*$ is a basis. 
	Thus, by replacing $\zz^*$ by an integer vector in $\poly(\overline{\AA}, \overline{\bb}) \setminus \{\zz^*\}$ finitely many times, we may assume that $\poly(\overline{\AA}, \overline{\bb}) \cap \Z^n = \{\zz^*\}$.
	\hfill\halmos
\endproof

\medskip

Observe that $\Delta^{\aalpha}\left(\AA\right) = \Delta^{\aalpha}(\overline{\AA})$.

\medskip

\proof{Proof of Theorem~\ref{thm:main_thm_linear_functional}.}
Suppose $\xx^*$ is an optimal vertex of $\LP(\AA, \bb, \cc)$. We apply Lemma \ref{lem:only_one_integer_vector}. Note that our bounds on $\prox_d(\AA,\bb,\aalpha)$ do not depend on the constraint matrix and right hand side. So we assume without loss of generality that $\poly\left(\AA,\bb\right)\cap \Z^n = \lbrace \zz^*\rbrace$.
%
%
Translating the instance, we may further assume that $\zz^*=\zero$, so that our objective is now to show $\bigl|\aalpha^\top\xx^{*}\bigr|<\frac{4n+2}{9}\cdot\Delta^{\objfn}\left(\AA\right)$.
%
%

Recall, that 
\[
\max_{\xx\in\poly}\aalpha^\top\xx=\prox_d(\aalpha)\cdot\Delta^{\aalpha}\left(\AA\right).
\]
By Lemma~\ref{lem:prox_dim_1_2_3}, $\prox_1(\aalpha)<1$, and since $n\geq 2$ we may assume $d\geq2$. We write $d=3a+2b$, where $a,b$ are non-negative integers, and we further specify
\[ 
a=\frac{d}{3}-2\cdot\left\{ -\frac{d}{3}\right\} \quad\text{and}\quad b=3\cdot\left\{ -\frac{d}{3}\right\}.
\]
where $\{x\}:=x-\lfloor x\rfloor$ denotes the fractional part of $x\in\R$. Applying Lemma~\ref{lem:proximity_template}, then Lemma~\ref{lem:prox_dim_1_2_3}, then the fact $d\leq n$, we get
\[
\prox_d(\aalpha)\leq\prox_3(\aalpha)\cdot a+\prox_2(\aalpha)\cdot b < \frac{4}{9}\cdot d + \frac{1}{3}\cdot\left\{ -\frac{d}{3}\right\} \leq \frac{4d+2}{9}\leq \frac{4n+2}{9}
\]
which implies
\begin{equation}
	\label{ineq:general:proof:main:thm}
	\max_{\xx\in\poly}\aalpha^\top\xx=\prox_d(\aalpha)\cdot\Delta^{\aalpha}\left(\AA\right)< \frac{4n+2}{9}\cdot\Delta^{\aalpha}\left(\AA\right).
\end{equation}
\hfill\halmos
\endproof

The proof of Theorem \ref{thm:main_thm_n_2} follows directly from inequality (\ref{ineq:general:proof:main:thm}) by setting $\aalpha=s\ee_i$ for $s\in\{-1,1\}$ and $i\in\{1,2,\ldots,n\}$. Note that under these assumptions $\Delta^{\aalpha}\left(\AA\right) \leq \Delta_{n-1}(\AA)$.

\section{Proof of Theorem \ref{thm:facet_width}.}
This section is devoted to outline a construction that allows us to derive Theorem \ref{thm:facet_width} from Theorem \ref{thm:main_thm_linear_functional}. In order to make this link precise, we define for a fixed full-column-rank matrix $\AA\in\Z^{m \times n}$ the parameter
\begin{align*}
		\pi\left(\AA\right) := \max_{\substack{\bb\in\Z^m \\ \text{ s.t. } \\ \poly\left(\AA,\bb\right)\cap\Z^n \neq \emptyset}} \max_{\substack{\xx^*\in \poly\left(\AA,\bb\right) \\ \text{ vertex}}} ~~\min_{\zz^*\in \poly\left(\AA,\bb\right) \cap \Z^n} ~\norm{\AA\left(\xx^* - \zz^*\right)}_\infty.
\end{align*}

The key connection between the lattice width and $\pi\left(\AA\right)$ is highlighted in the statement below.

\begin{lemma}\label{lem:link_facet_width_gap}
	Let $n\geq2$ and $\bb\in\Z^m$ such that $\poly\left(\AA,\bb\right)$ is a full-dimensional lattice-free polyhedron and each row of $\AA$ is facet-defining. Then, there exists a row $\aa$ of $\AA$ such that
	\begin{align*}
		w^\mathcal{\aa}\left(\poly\left(\AA,\bb\right)\right)\leq \pi\left(\AA\right) - 1.
	\end{align*}
\end{lemma}
\proof{Proof.}
	Throughout the proof, we abbreviate $\ee_J := \sum_{j\in J}\ee_j$ for $J\subseteq\lbrack m\rbrack$. Let $k\in \mathbb{N}$ be given such that $\poly\left(\AA,\bb + k\one\right)$ contains integer vectors.
	We define for each $\zz \in \poly\left(\AA,\bb + k\one\right)\cap \Z^n$ the set $I(\zz) := \lbrace i\in \lbrack m\rbrack : \aa_i^\top\zz\geq \bb_i + 1\rbrace$.
	
	We choose a $\overline{\zz}\in \poly\left(\AA,\bb + k\one\right)\cap \Z^n$ such that $|I(\zz)|$ is minimal among all integer vectors in $\poly\left(\AA,\bb + k\one\right)$. For sake of brevity, we set $I:= I(\overline{\zz})$.
	
	In the following, we analyze $\poly\left(\AA,\bb + k \ee_I\right)$. This polyhedron is not lattice-free since $\overline{\zz}\in \poly\left(\AA,\bb + k \ee_I\right)$. Furthermore, the minimality of $|I|$ implies
	\begin{align}
		\label{proof_min_key_property}
		\AA_{I}\zz\geq \bb_I + \one
	\end{align}
	for all $\zz\in \poly\left(\AA,\bb + k\ee_I\right)\cap\Z^n$. 
	
	Pick $i\in I$, observe that $I\neq \emptyset$ as $\poly\left(\AA,\bb\right)$ is lattice-free. Choose a vertex $\tilde{\xx}$ which minimizes $\aa_i^\top\xx$ over $\poly\left(\AA,\bb\right)$ and a vertex $\xx^*$ which minimizes $\aa_i^\top\xx$ over $\poly\left(\AA,\bb + k \ee_I\right)$. Since $\poly(\AA,\bb)$ is not necessarily bounded, it is not obvious why these vertices exist in the first place. We claim that our choice of $\aa_i$ implies that: If $\min\aa_i^\top\xx$ is unbounded over $\poly\left(\AA,\bb + k \ee_I\right)$, then there exists some $\rr\in\Z^n$ with $\AA\rr\leq\zero$ such that $\aa_i^\top\rr \leq -1$. This yields $\aa_i^\top(\zz-\lambda\rr)\leq \bb_i$ for some $\zz\in \poly\left(\AA,\bb + k\ee_I\right)\cap\Z^n$ and some large enough $\lambda\in \Z_{\geq 0}$, contradicting (\ref{proof_min_key_property}). So we have $\min \aa_i^T\xx$ over $\poly(\AA,\bb + k\ee_I)$ is bounded which also implies boundedness over $\poly(\AA,\bb)$ as $\poly(\AA,\bb)\subseteq \poly(\AA,\bb + k\ee_I)$. 
	
	There exists $\zz^*\in \poly\left(\AA,\bb + k \ee_I\right)\cap\Z^n$ such that 
	\begin{align}
		\label{proof_min_prox}
		\aa_i^\top(\zz^*-\xx^*)\leq \pi\left(\AA\right).
	\end{align}
	Let $\overline{\yy} \in \poly\left(\AA,\bb\right)$ be a vertex maximizing $\aa_i^\top\xx$. So we have $w^{\aa_i}\left(\poly\left(\AA,\bb\right)\right)= \aa_i^\top(\tilde{\yy} - \tilde{\xx})$. We obtain
	
	\begin{align*}
		w^{\aa_i}\left(\poly\left(\AA,\bb\right)\right) = \aa_i^\top(\tilde{\yy} - \tilde{\xx})
		\leq \bb_i - \aa_i^\top\xx^* 
		\leq \aa_i^\top\zz^* - \aa_i^\top\xx^* - 1
		\leq \pi\left(\AA\right) - 1,
	\end{align*}
	where the first inequality comes from the fact that $\poly\left(\AA,\bb\right)\subseteq \poly\left(\AA,\bb + k\ee_I\right)$. We use (\ref{proof_min_key_property}) for the second inequality and (\ref{proof_min_prox}) for the third inequality.
	\hfill\halmos
\endproof

\medskip

\proof{Proof of Theorem~\ref{thm:facet_width}.}
	Our strategy is to apply Lemma \ref{lem:link_facet_width_gap}. We can bound $\pi\left(\AA\right)$ using Theorem \ref{thm:main_thm_linear_functional}. 
	In order to do so, let $\pi(\AA)$ be attained for $\poly(\AA,\hat{\bb})$ with $\hat{\bb}\in\Z^m$, a vertex $\xx^*\in\poly(\AA,\hat{\bb})$, and some integral vector in $\poly(\AA,\hat{\bb})$. Thus, $\pi(\AA)\leq \norm{\AA\left( \xx^* - \zz^*\right)}_\infty$ for all $\zz^*\in \poly(\AA,\hat{\bb})\cap\Z^n$.
	
	Since $\xx^*$ is a vertex, there exists some $\cc\in\Q^n$ such that $\cc^\top\xx$ is maximized at $\xx^*$ over $\poly(\AA,\hat{\bb})$. This enables us to apply Lemma \ref{lem:only_one_integer_vector} and obtain a polyhedron $\poly(\overline{\AA} ,\overline{\bb})\subseteq \poly(\AA,\hat{\bb})$ with $\poly(\overline{\AA},\overline{\bb})\cap\Z^n = \lbrace \zz^*\rbrace$ and $\xx^*\in\poly(\overline{\AA},\overline{\bb})$, where $\overline{\bb}$ is integral and the rows of $\overline{\AA}$ are a subset of the rows of $\AA$ and their negatives. We analyze the distance between $\xx^*$ and $\zz^*$ with respect to $\poly(\overline{\AA} ,\overline{\bb})$ via Theorem \ref{thm:main_thm_linear_functional}. Observe that $\zz^*$ is always the closest integer vector in $\poly(\overline{\AA},\overline{\bb})$ to $\xx^*$ with respect to some arbitrary $\aalpha$ since $\zz^*$ is the only integer vector in $\poly(\overline{\AA},\overline{\bb})$.
	
	Choose $\aalpha$ in the statement of Theorem \ref{thm:main_thm_linear_functional} to be a row of $\AA$, say $\aa_{i}$ for some $i\in\lbrack m\rbrack$. This leads to
	\begin{align*}
		\left| \aa_{i}^\top(\xx^* - \zz^*)\right| < \frac{4n+2}{9}\cdot\Delta_n(\AA).
	\end{align*}
	Note that $\Delta^{\aa_{i}}(\AA) \leq \Delta_n(\AA)$ and $\Delta_n(\AA) = \Delta_n(\overline{\AA})$. This procedure works for every row of $\AA$. Hence, we conclude $\norm{\AA\left( \xx^* - \zz^*\right)}_\infty < \frac{4n+2}{9}\cdot\Delta_n(\AA)$. Therefore, we get 
	\begin{align*}
		\pi\left(\AA\right)\leq \norm{\AA\left( \xx^* - \zz^*\right)}_\infty < \frac{4n+2}{9}\cdot\Delta_n(\AA)	
	\end{align*}
	and the claim follows from Lemma \ref{lem:link_facet_width_gap}.\hfill\halmos
\endproof

\section{The $\lbrace 0,\pm k, \pm 2k\rbrace$-case.}\label{secDimensionFree}
In this section, we prove Theorem \ref{thm:dim_free_bounds}, that is, bounds on the proximity and facet width of lattice-free polyhedra which are independent of the dimension.

For this purpose, we need the notion of Graver bases. Given a full-column-rank matrix $\AA\in\Z^{m\times n}$, we define the cone $\ccone\left(\AA\right) := \lbrace \xx\in\R^n : \AA\xx\geq \zero\rbrace$. There exists a unique minimal set $\mathcal{H}\left(\AA\right)\subseteq \ccone(\AA)\cap\Z^n$ such that every element in $\ccone\left(\AA\right)\cap \Z^n$ is a non-negative integral combination of the elements in $\mathcal{H}\left(\AA\right)$. This set is called \emph{Hilbert basis of $\ccone(\AA)$} and its elements are referred to as \emph{Hilbert basis elements}. Then, the \emph{Graver basis of $\AA$} is given by 
\begin{align*}
	\mathcal{G}\left(\AA\right) := \bigcup_{\SSv\in D} \mathcal{H}\left(\SSv\AA\right),
\end{align*}
where $D$ is the set of all diagonal $m\times m$-matrices with $\pm 1$ entries on the diagonal. Note that $\ccone\left(\AA\right) = \lbrace \zero\rbrace$ implies $\mathcal{H}\left(\AA\right) = \emptyset$. We refer to the elements of $\mathcal{G}\left(\AA\right)$ as \emph{Graver basis elements}. 

The Hilbert basis elements satisfy an important property: They are precisely the irreducible elements in $\ccone\left(\AA\right)$, i.e., given $\yy^1,\yy^2\in \ccone\left(\AA\right)\cap \Z^n$ with $\hh = \yy^1 + \yy^2$ for $\hh\in\mathcal{H}\left(\AA\right)$, we have that either $\yy^1 = \zero$ or $\yy^2 = \zero$. This is the case if and only if $\spindle\left(\AA,\hh\right) \cap \Z^n = \lbrace \zero, \hh\rbrace$. 

The main result is based on taking suitable Graver basis steps in a certain polytope. Since we aim to measure the length of these steps with respect to some $\aalpha\in \Z^n$, we define $\tilde{\kappa}_n\left(\AA,\aalpha\right)$ to be the minimum number such that 
\begin{align*}
	\left|\aalpha^\top\gg\right|\leq \tilde{\kappa}_n\left(\AA,\aalpha\right)\cdot \Delta^{\aalpha}\left(\AA\right)
\end{align*}
for all $\gg\in \mathcal{G}\left(\AA\right)$.

Note that in the following we work with polyhedra $\poly(\AA\BB,\bb)$, where $\BB\in\Z^{n\times n}$ is invertible. In order to highlight the dependence on $\AA$ and $\BB$, we write $\kappa_n\left(\AA,\bb,\aalpha\right)$ but allow for rows of $-\AA$ in the definition of $\kappa_n\left(\AA,\bb,\aalpha\right)$.
\begin{lemma}\label{lem:bound_gap_transformation}
	Let $\aalpha\in\Z^n\backslash\lbrace \zero\rbrace$, $\AA\in\Z^{m\times n}$ have full column rank, $\BB\in\Z^{n\times n}$ be invertible, and $\bb\in\Z^m$ such that $\poly\left(\AA\BB,\bb\right)\cap\Z^n = \lbrace \zero \rbrace$. Let $\xx^*$ be a vertex of $\poly\left(\AA\BB,\bb\right)\cap\Z^n$. Then,
	\begin{align*}
		\left|\aalpha^\top\xx^*\right| \leq \frac{\kappa_n\left(\AA,\overline{\bb},\aalpha\right) + \tilde{\kappa}_n\left(\AA,\aalpha\right)(\left|\det\BB\right| - 1)}{\left|\det\BB\right|}\cdot\Delta^{\aalpha}\left(\AA\BB\right)
	\end{align*}
	for some integral vector $\overline{\bb}$.	
\end{lemma}
\proof{Proof.}
	Observe that $\BB\cdot \poly\left(\AA\BB,\bb\right) = \poly\left(\AA,\bb\right)$ and define the lattice $\Lambda : = \BB\Z^n$. Note that $\Lambda\subseteq \Z^n$ and $\poly\left(\AA,\bb\right)\cap\Lambda = \lbrace \zero\rbrace$. 
	Moreover, we set $\bbeta := \left|\det\BB\right|\BB^{-\top}\aalpha\in\Z^n$ and $\yy^* := \BB\xx^*$. Our aim is to bound 
	\begin{align*}
		\left|\aalpha^\top\xx^*\right| = \frac{1}{\left|\det\BB\right|}\left|\bbeta^\top\yy^*\right|.
	\end{align*}

	Since $\xx^*$ is a vertex of $\poly\left(\AA\BB,\bb\right)$, $\yy^* = \BB\xx^*$ is a vertex of $\poly\left(\AA,\bb\right)$. 
	So there exists some $\cc\in\Q^n$ which is maximized by $\yy^*$ over $\poly\left(\AA,\bb\right)$. 
	We know that ILP($\AA,\bb,\cc$) is feasible as $\zero\in \poly\left(\AA,\bb\right)\cap\Z^n$. 
	Thus, we can apply Lemma \ref{lem:only_one_integer_vector}: 
	there exists $\zz^{(0)}\in \poly\left(\AA,\bb\right)\cap\Z^n$, an optimal solution of ILP($\AA,\bb,\cc$), and $\poly(\overline{\AA},\overline{\bb})$ with $\poly(\overline{\AA},\overline{\bb})\cap\Z^n = \lbrace \zz^{(0)}\rbrace$, where the rows of $\overline{\AA}$ are rows of $\AA$ or their negatives. 
	Since the definition of $\kappa_n\left(\AA,\bb,\aalpha\right)$ is translation-invariant, we conclude
	\begin{align}\label{ineq_bound_gap_transformation_gap}
		\left|\bbeta^\top(\yy^*-\zz^{(0)})\right|\leq \kappa_n\left(\AA,\overline{\bb},\aalpha\right)\cdot\Delta^{\bbeta}\left(\AA\right).
	\end{align}
	If $\zz^{(0)} = \zero$, then we are done and the claim follows from $\Delta^{\bbeta}\left(\AA\right) = \Delta^{\aalpha}\left(\AA\BB\right)$. Suppose that $\zz^{(0)}\neq\zero$.
	We go back to analyzing $\poly\left(\AA,\bb\right)$ and pass to the spindle $\spindle\left(\AA,\zz^{(0)}\right)\subseteq \poly\left(\AA,\bb\right)$. We claim that there exists $\zz^{(1)}\in \spindle\left(\AA,\zz^{(0)}\right)\cap\Z^n$ such that $\zz^{(1)}\neq\zz^{(0)}$ and $\spindle\left(\AA,\zz^{(0)}-\zz^{(1)}\right)\cap\Z^n = \lbrace\zero,\zz^{(0)} - \zz^{(1)}\rbrace$. 
	Recall that $D$ denotes the set of all $m\times m$ diagonal matrices with $\pm1$ entries on the diagonal. 
	Let $\SSv\in D$ be such that $\spindle\left(\AA,\zz^{(0)}\right) = \lbrace \xx \in\R^n : \zero\leq \SSv\AA\xx\leq\SSv\AA\zz^{(0)}\rbrace$. Further, let $\zz\in \spindle\left(\AA,\zz^{(0)}\right)\cap\Z^n\backslash\lbrace \zz^{(0)}\rbrace$. If $\spindle\left(\AA,\zz^{(0)}-\zz\right)\cap\Z^n = \lbrace \zero,\zz^{(0)}-\zz\rbrace$, we set $\zz^{(1)} := \zz$. Otherwise, there exists $\tilde{\zz}\in\spindle\left(\AA,\zz^{(0)}\right)\cap\Z^n\backslash\lbrace \zz^{(0)},\zz\rbrace$ such that $\SSv\AA\zz\leq \SSv\AA\tilde{\zz}$. We pass to $\spindle\left(\AA,\zz^{(0)}-\tilde{\zz}\right)$ and iterate this procedure. Note that $\left|\spindle\left(\AA,\zz^{(0)}-\zz\right)\cap\Z^n\right| > \left|\spindle\left(\AA,\zz^{(0)}-\tilde{\zz}\right)\cap\Z^n\right|$ as $\zz^{(0)}-\zz\notin\spindle\left(\AA,\zz^{(0)}-\tilde{\zz}\right)\cap\Z^n$. Hence, our procedure terminates with some $\zz^{(1)}\in \spindle\left(\AA,\zz^{(0)}\right)\cap\Z^n\backslash\lbrace \zz^{(0)}\rbrace$ such that $\spindle\left(\AA,\zz^{(0)}-\zz^{(1)}\right)\cap\Z^n = \lbrace \zero,\zz^{(1)}\rbrace$.
	
	This choice of $\zz^{(1)}$ guarantees that $\zz^{(0)}-\zz^{(1)}\in \ccone\left(\SSv\AA\right)\cap\Z^n$ is irreducible and, thus, $\zz^{(0)}-\zz^{(1)}\in\mathcal{H}\left(\SSv\AA \right)$. Hence, we have $\zz^{(0)}-\zz^{(1)}\in\mathcal{G}\left(\AA\right)$ by the definition of Graver bases. Therefore, we get
	\begin{align}\label{ineq_bound_gap_transformation_graver}
		\left|\bbeta^\top(\zz^{(0)}-\zz^{(1)})\right|\leq \tilde{\kappa}_n\left(\AA,\aalpha\right)\cdot\Delta^{\bbeta}\left(\AA\right).
	\end{align}
	We pass to $\spindle\left(\AA,\zz^{(1)}\right)\subseteq \spindle\left(\AA,\zz^{(0)}\right)$ and repeat the procedure until $\zz^{(s)} =\zero$ for some integer $s\geq 1$. 
	
	We claim that $s\leq\left|\det\BB\right| - 1$. If $s\geq \left|\det\BB\right|$, then the pigeonhole principle gives us that either there exists a vector $\zz^{(i)}\in\Lambda$ for some $i\in\lbrace 0,\ldots,s - 1\rbrace$ or there exist vectors such that $\zz^{(k)}-\zz^{(l)}\in\Lambda$ for $k,l\in\lbrace 0,\ldots,s - 1\rbrace$ with $l>k$. The first case is not possible since $\zero\neq\zz^{(i)}\in \spindle\left(\AA,\zz^{(0)}\right)\cap\Lambda\subseteq \poly\left(\AA,\bb\right)\cap\Lambda$, contradicting $\poly\left(\AA,\bb\right)\cap\Lambda = \lbrace\zero\rbrace$. The second case leads to the same contradiction as $\zero\neq\zz^{(k)} - \zz^{(l)}\in \spindle\left(\AA,\zz^{(0)}\right)\cap\Lambda$. Thus, we have $s\leq \left|\det\BB\right| - 1$.
	
	As a result, we obtain
	\begin{align*}
		\left|\aalpha^\top\xx^*\right| &= \frac{1}{\left|\det\BB\right|}\left|\bbeta^\top\yy^*\right| \\
		&\leq \frac{1}{\left|\det\BB\right|}\left(\left|\bbeta^\top(\yy^*-\zz^{(0)})\right| + \sum_{i = 0}^{s - 1}\left|\bbeta^\top(\zz^{(i)}-\zz^{(i+1)})\right|\right) \\
		&\leq\frac{\kappa_n\left(\AA,\overline{\bb},\aalpha\right) + \tilde{\kappa}_n\left(\AA,\aalpha\right)(\left|\det\BB\right| - 1)}{\left|\det\BB\right|}\cdot\Delta^{\aalpha}\left(\AA\BB\right),
	\end{align*}
	where we use (\ref{ineq_bound_gap_transformation_gap}), (\ref{ineq_bound_gap_transformation_graver}), and $s\leq \left|\det\BB\right| - 1$ for the last inequality. Moreover, we exploit that $\Delta^{\bbeta}\left(\AA\right) = \Delta^{\aalpha}\left(\AA\BB\right)$.
\hfill\halmos
\endproof

\medskip

We want to utilize Lemma \ref{lem:bound_gap_transformation}. 
Therefore, we need upper bounds on $\kappa_n(\AA,\overline{\bb},\aalpha)$ and $\tilde{\kappa}_n(\AA,\aalpha)$. If $\AA$ is unimodular, we have $\kappa_n(\AA,\overline{\bb},\aalpha) = 0$ as every vertex of $\poly(\AA,\overline{\bb})$ is integral. 
In order to determine $\tilde{\kappa}_n\left(\AA,\aalpha\right)$ when $\AA$ is unimodular, we claim that  $\left|\aalpha^\top\rr\right|\leq\Delta^{\aalpha}\left(\AA\right)$, where $\rr$ is a primitive vector on an extreme ray of $\ccone\left(\AA\right)$, that is, $\rr\in\Z^n$ and $\gcd\rr = 1$. More specifically, there exists $I\subseteq \lbrack m\rbrack$ with $|I|=n-1$ such that $\text{rank}~\AA_{I} = n -1$ and, up to a sign,
\begin{align*}
	\rr_i = \frac{1}{\gcd\AA_I}(-1)^i\det\AA_{I,\lbrack n\rbrack \backslash i}
\end{align*}
for all $i\in\lbrack n\rbrack$, where $\AA_{I,J}$ denotes the matrix with rows indexed by $I$ and columns indexed $J$ for $I\subseteq\lbrack m\rbrack$ and $J\subseteq\lbrack n\rbrack$. So
\begin{align}\label{ineq:dim_free_primitive}
	\left|\aalpha^\top\rr\right|\leq\Delta^{\aalpha}\left(\AA\right)
\end{align}
by Laplace expansion. 

It is well-known that the Hilbert basis elements coincide with the primitive vectors on the extreme rays of $\ccone\left(\AA\right)$ if $\AA$ is unimodular; see, e.g., \cite[Proposition 8.1]{Sturmfels1995GrobnerBA}. Thus, we have $\left|\aalpha^\top\hh\right|\leq\Delta^{\aalpha}\left(\AA\right)$ for all $\hh \in \mathcal{H}\left(\AA\right)$ by (\ref{ineq:dim_free_primitive}). This extends naturally to the Graver basis and we conclude $\tilde{\kappa}_n\left(\AA,\aalpha\right)\leq 1$.

In order to prove Theorem \ref{thm:dim_free_bounds}, we still need bounds on $\kappa_n\left(\AA,\overline{\bb},\aalpha\right)$ and $\tilde{\kappa}_n\left(\AA,\aalpha\right)$ when $\AA$ is \emph{bimodular}, i.e., $\Delta_n\left(\AA\right) = 2$.

\begin{lemma}\label{lem:dim_free_bimodular}
	Let $\aalpha\in\Z^n\backslash\lbrace \zero\rbrace$, $\AA\in\Z^{m\times n}$ be bimodular, and $\bb\in\Z^m$ such that $\poly\left(\AA,\bb\right)\cap\Z^n = \lbrace \zero \rbrace$. Then we have
	\begin{enumerate}
		\item $\kappa_n\left(\AA,\bb,\aalpha\right)\leq \frac{1}{2}$ and
		\item $\tilde{\kappa}_n\left(\AA,\aalpha\right)\leq 1$.
	\end{enumerate}
\end{lemma}
\proof{Proof.}
	We begin with $\kappa_n\left(\AA,\bb,\aalpha\right)$. Let $\xx^*$ be an arbitrary vertex of $\poly\left(\AA,\bb\right)$. If $\xx^*\in\Z^n$, then $\xx^*=\zero$ because $\poly\left(\AA,\bb\right)\cap\Z^n = \lbrace \zero \rbrace$ and, in particular, $\left|\aalpha^\top\xx^*\right|= 0$. 
	
	So suppose that $\xx^*\notin\Z^n$. Additionally, assume that without loss of generality $\dim \poly\left(\AA,\bb\right) = n$ by Lemma \ref{lem:lifting_lemma}. This assumption allows us to apply a result by Chirkov and Veselov \cite[Theorem 2]{VC2009}: There exists $\zz\in \poly\left(\AA,\bb\right)\cap\Z^n$ such that $\zz$ and $\xx^*$ lie on an edge of $\poly\left(\AA,\bb\right)$. Since $\poly\left(\AA,\bb\right)\cap\Z^n=\lbrace\zero\rbrace$, we have $\zz = \zero$. So $\xx^*\in\ker\AA_{I}$ for some $I\subseteq\lbrack m\rbrack$ with $|I| = n - 1$ and $\text{rank}~\AA_{I} = n - 1$. Since the line segment $\lbrack \zero, \xx^*\rbrack$ is an edge of $\poly\left(\AA,\bb\right)$ and $\poly\left(\AA,\bb\right)\cap\Z^n = \lbrace \zero\rbrace$, there is no non-zero integer vector contained in $\lbrack \zero, \xx^*\rbrack$. Additionally, we have $2\xx^*\in\ker\AA_{I}\cap\Z^n$ by Cramer's rule. We conclude that $2\xx^*$ is primitive and get
	\begin{align*}
		\left|2\aalpha^\top\xx^*\right|\leq \Delta^{\aalpha}\left(\AA\right)
	\end{align*}
	by (\ref{ineq:dim_free_primitive}). Dividing by two yields the first part of the statement.
	
	For the second statement, we utilize a structural result about the Hilbert basis of cones defined by bimodular matrices \cite[Theorem 1.5]{HKW2022}: Every $\hh\in\mathcal{H}\left(\AA\right)$ either is a primitive vector on an extreme ray or can be expressed as $\frac{1}{2}\rr^1+\frac{1}{2}\rr^2$, where $\rr^1$ and $\rr^2$ are integral vectors on extreme rays such that $\left|\aalpha^\top\rr^i\right|\leq\Delta^{\aalpha}\left(\AA\right)$ for $i=1,2$. In first case, we follow that $\left|\aalpha^\top\hh\right|\leq \Delta^{\aalpha}\left(\AA\right)$ by (\ref{ineq:dim_free_primitive}) and in the second case we draw the same conclusion from the triangle inequality. This holds for every bimodular cone and, therefore, generalizes naturally to $\tilde{\kappa}_n\left(\AA,\aalpha\right)\leq 1$.
	\hfill\halmos
\endproof

\medskip
We are in the position to prove our main result.

\medskip
\proof{Proof of Theorem \ref{thm:dim_free_bounds}.}
	We apply Lemma \ref{lem:only_one_integer_vector}. Note that this does not alter the property that every $n\times n$ minor is contained in $\lbrace 0,\pm k,\pm 2k\rbrace$. So we assume without loss of generality that $\poly\left(\AA,\bb\right)\cap\Z^n = \lbrace \zero\rbrace$. 
	
	We can decompose $\AA = \TT\BB$ such that $\left|\det\BB\right| = k$ and the $n\times n$ minors of $\TT$ are contained in $\lbrace 0,\pm 1,\pm 2\rbrace$. Thus, $\TT$ is either unimodular or bimodular. 
	Our aim is to apply Lemma \ref{lem:bound_gap_transformation}. 
	Observe that $\tilde{\kappa}_n\left(\TT,\aalpha\right)\leq 1$ by the previous discussion and Lemma \ref{lem:dim_free_bimodular}. 
	Recall that $\kappa_n\left(\TT,\overline{\bb},\aalpha\right) = 0$ when $\TT$ is unimodular and $\kappa_n\left(\TT,\overline{\bb},\aalpha\right) \leq \frac{1}{2}$ if $\TT$ is bimodular by Lemma \ref{lem:dim_free_bimodular}. 
	Together we have $\kappa_n\left(\TT,\overline{\bb},\aalpha\right) \leq \frac{\Delta_n\left(\TT\right) - 1}{\Delta_n\left(\TT\right)}$. 
	By Lemma \ref{lem:bound_gap_transformation}, we obtain
	\begin{align}\label{ineq:bound_dim_free}
		\left|\aalpha^\top\xx^*\right|\nonumber&\leq \frac{\kappa_n\left(\TT,\overline{\bb},\aalpha\right) + \tilde{\kappa}_n\left(\TT,\aalpha\right)(\left|\det\BB\right| - 1)}{\left|\det\BB\right|}\cdot\Delta^{\aalpha}\left(\AA\right)\\\nonumber
		&\leq \left(1 - \frac{1}{\Delta_n\left(\TT\right)\left|\det\BB\right|}\right)\cdot \Delta^{\aalpha}\left(\AA\right)\\ 
		&= \frac{\Delta_n\left(\AA\right) - 1}{\Delta_n\left(\AA\right)}\cdot \Delta^{\aalpha}\left(\AA\right)
	\end{align}
	as $\Delta_n\left(\TT\right)\left|\det\BB\right| = \Delta_n\left(\TT\BB\right) = \Delta_n\left(\AA\right)$.
	
	The proximity bounds of the statement follow directly by setting $\aalpha = \ee_i$ for $i\in\lbrace 1,2,\ldots,n\rbrace$ and the observation that $\Delta^{\aalpha}\left(\AA\right)\leq\Delta_{n-1}\left(\AA\right)$ for this choice of $\aalpha$.
	
	In order to bound the facet width, we follow the same strategy as in the proof of Theorem \ref{thm:facet_width}. In particular, we want to use Lemma \ref{lem:link_facet_width_gap}. 
	So our aim is to bound $\pi\left(\AA\right)$. 
	Applying the reduction from the proof of Theorem \ref{thm:facet_width}, we assume that $\zz^*$ is the only integer vector contained in our polyhedron. 
	We choose $\aalpha$ to be a row of $\AA$, say $\aa_i$ for some $i\in\lbrack m\rbrack$. 
	By (\ref{ineq:bound_dim_free}), we get
	\begin{align*}
		\left|\aa_i^\top\left(\xx^*-\zz^*\right)\right|\leq \Delta_n\left(\AA\right) - 1,
	\end{align*}
	where we use that $\Delta^{\aa_i}\left(\AA\right)\leq\Delta_n\left(\AA\right)$. 
	This bound holds for all $i\in\lbrack m\rbrack$ and hence $\Vert \AA\left(\xx^* - \zz^*\right)\Vert_\infty\leq \Delta_n\left(\AA\right) - 1$. 
	Therefore, we get $\pi(\AA)\leq \Delta_n\left(\AA\right) - 1$ and the claim follows from Lemma \ref{lem:link_facet_width_gap}.
	\hfill\halmos
\endproof

\begin{APPENDICES}

\section{The area of a polygon containing its rotated polar.}\label{appAreaPolygonPolar}

In this appendix we prove Lemma~\ref{lem:polar_area_lower_bound}, which states that any polygon $\qoly\subseteq\R^2$ satisfying $\rot\qoly\subseteq\qoly^{\circ}$ has $\vol_2(\qoly^{\circ})\geq3$. 
Recall that $\rot:\R^{2}\rightarrow\R^{2}$ denotes the $90^{\circ}$ counterclockwise
rotation in $\R^{2}$. 
We frequently use here the fact that for all $\xx\in\R^{2}$ we have $(\rot\vv)^\top\xx=\det(\vv,\xx)$.

It turns out that we only need to consider closed convex polygons $\poly\subseteq\R^{2}$ satisfying
the equality $\rot\poly=\poly^{\circ}$. 
Note that any such $\poly$ must contain the origin. 
Such polygons do exist, for example  suitably scaled regular ($4k+2$)-gons for $k\geq1$. 
More general examples which are not polygons include lines through the origin and the unit Euclidean ball. 

A thorough analysis of polytopes which are linearly equivalent to their polars was undertaken in~\cite{J2021}. 
We use many of the results from this work, in particular the modification techniques from~\cite[Section 7]{J2021}, and provide corresponding citations where appropriate. 
These techniques do not directly apply to our setting, as they concern polytopes satisfying $-\poly=\poly^{\circ}$. 
Thus, the proofs here are self-contained. 
However, this difference seems to be superficial at first glance, and it would be interesting to unite the two points of view.

We use the following notation in our proof: letters $\poly$ and $\qoly$ refer to convex sets in the plane, typically polygons. Greek letters such as $\rot,\gamma$ refer to linear transformations in the plane. 
Letters $\uu,\vv,\ww$ refer to points in the plane, typically vertices of polygons. 
For a set $S$ we write $\pm S$ as shorthand for $S\cup-S$. 
For a point $\vv$ in the plane we write $\pm\vv$ as shorthand for $\vv,-\vv$. 
Closed (resp. open) line segments in the plane with ends $\uu,\vv$ are denoted by $[\uu,\vv]$ (resp. $(\uu,\vv)$). 
\begin{proposition}[{cf. \cite[Theorem 3.3]{J2021}}]
\label{prop:centrally-symmetric}Suppose $\rot\poly=\poly^{\circ}$.
Then $\poly$ is centrally-symmetric.
\end{proposition}

\proof{Proof.}
For any non-singular linear transformation $\gamma$ on $\R^{2}$,
we have $(\gamma\poly)^{\circ}=\gamma^{-\top}\poly^{\circ}$. 
Since $\rot^{-\top}=\rot$, we have
\[
-\poly=\rot^{2}\poly=\rot(\rot\poly)=\rot\poly^{\circ}=(\rot^{-\top}\poly)^{\circ}=(\rot\poly)^{\circ}=(\poly^{\circ})^{\circ}=\poly.\tag*{\halmos}
\]
\endproof

\begin{proposition}
\label{prop:classification}Suppose $\rot\poly=\poly^{\circ}$. 
Then either $\poly$ is a line through the origin, or $\poly$ is bounded and full-dimensional.
\end{proposition}

\proof{Proof.}
Suppose $\poly$ is not a line through the origin. 
We rule out that $\poly$ is contained in a line through the origin. 
If this were the case, then by Proposition~\ref{prop:centrally-symmetric}, $\poly$ must be bounded, which implies $\poly^{\circ}$ is full-dimensional.
But this would contradict $\tau\poly=\poly^{\circ}$. 
So $\poly$ is not contained in a line through the origin. 
Since $\zero\in\poly$, $\poly$ must therefore be full-dimensional. 
If $\poly$ were unbounded, then $\poly^{\circ}$ would be lower-dimensional, but again this would contradict $\tau\poly=\poly^{\circ}$.\hfill\halmos
\endproof
\begin{proposition}
\label{prop:absdet1-invariant}Suppose $\rot\poly=\poly^{\circ}$.
If $\gamma$ is a $2\times2$ matrix with determinant $\pm1$, then
$\gamma\poly$ also has this property.
\end{proposition}

\proof{Proof.}
Let $s=\det(\gamma)$. One quickly verifies that $\gamma^\top\rot\gamma=s\rot$.
We have
\[
\rot(\gamma\poly)=s\gamma^{-\top}\rot\poly=s\gamma^{-\top}\poly^{\circ}=s(\gamma\poly)^{\circ}=(\gamma(s\poly))^{\circ}=(\gamma\poly)^{\circ},
\]
where the last equality holds by central symmetry.\hfill\halmos
\endproof
\begin{definition}
Let $\vv\in\R^{2}$ be non-empty. We define the line $L_{\vv}$, the
half-space $H_{\vv}$, and the strip $S_{\vv}$ as
\[
\begin{array}{rclcl}
L_{\vv} & :=&\{ \xx\in\R^{2}:(\rot\vv)^\top\xx&=&1\} \\
H_{\vv} & :=&\{ \xx\in\R^{2}:(\rot\vv)^\top\xx&\leq&1\} \\
S_{\vv} & :=&\{ \xx\in\R^{2}:\left|(\rot\vv)^\top\xx\right|&\leq&1\} .
\end{array}
\]
Note that $\uu\in L_{\vv}\Leftrightarrow\vv\in L_{-\uu}$ and similarly
$\uu\in H_{\vv}\Leftrightarrow\vv\in H_{-\uu}$. On the other hand,
we have $\uu\in S_{\vv}\Leftrightarrow\vv\in S_{\uu}$.
\end{definition}

In the remainder of this manuscript we focus only on the case when
$\poly$ is a polygon. Let $\verts(\poly)$ and $\edges(\poly)$ denote
the set of vertices and edges of $\poly$, respectively.
\begin{proposition}[{cf.~\cite[Theorem~3.2]{J2021}}]
\label{prop:edge_vertex_bijection}Suppose $\rot\poly=\poly^{\circ}$
is a polygon. Then the map $\left[\uu,\vv\right]\mapsto L_{\uu}\cap L_{\vv}$
is a bijection from $\edges(\poly)$ to $\verts(\poly)$. The inverse
of this bijection is given by $\vv\mapsto\poly\cap L_{-\vv}$.
\end{proposition}

\proof{Proof.}
For $\vv\in\R^{2}$, define $\ell_{\vv}$ to be the line $\vv^\top\xx=1$.
The map $\left[\uu,\vv\right]\mapsto\ell_{\uu}\cap\ell_{\vv}$ is
a bijection from $\edges(\poly)$ to $\verts(\poly^{\circ})=\verts(\rot\poly)$,
with inverse $\vv\mapsto\poly\cap\ell_{\vv}$. 
The map $\vv\mapsto\rot\vv$ is a bijection from $\verts(\rot\poly)$ to $\verts(\poly)$ by central symmetry, with inverse $\vv\mapsto-\rot\vv$. 
Composing these two maps yields the desired bijection.
\hfill\halmos
\endproof
\begin{proposition}
\label{prop:at_least_six_vertices}Suppose $\rot\poly=\poly^{\circ}$
is a polygon. Then $\left|\verts(\poly)\right|\geq6$. 
\end{proposition}

\proof{Proof.}
Since $\poly$ is centrally symmetric, $\left|\verts(\poly)\right|$
is even, and hence $\left|\verts(\poly)\right|\geq4$. But we cannot
have $\left|\verts(\poly)\right|=4$. Suppose this were the case.
Then $\poly$ is a parallelogram by central symmetry. By Proposition~\ref{prop:absdet1-invariant},
we may apply a suitable determinant 1 linear transformation so that
$\poly$, and hence $\rot\poly$, is an origin-symmetric axis-aligned
square. But then $\poly^{\circ}$ is a two-dimensional cross-polytope.
Thus $\rot\poly\neq\poly^{\circ}$. So $\left|\verts(\poly)\right|\geq6$.\hfill\halmos
\endproof
\begin{proposition}
\label{prop:sufficient_condition_for_six_verts}Suppose $\rot\poly=\poly^{\circ}$
is a polygon, and there exists three consecutive vertices $\uu<\vv<\ww$
in the counterclockwise order such that $L_{\uu}$ contains both $\vv$
and $\ww$. Then $\left|\verts(\poly)\right|=6$.
\end{proposition}

\proof{Proof.}
Since $L_{\uu}$ contains both $\vv$ and $\ww$, we have $\poly\cap L_{\uu}=\left[\vv,\ww\right]$,
and hence $\uu$ is the intersection of the two edges $\poly\cap L_{-\vv}$
and $\poly\cap L_{-\ww}$ by Proposition~\ref{prop:edge_vertex_bijection}.
It follows that either $\left[\uu,\vv\right]=\poly\cap L_{-\vv}$
or $\left[\uu,\vv\right]=\poly\cap L_{-\ww}$. Since $\vv\notin L_{-\vv}$,
we must have $\left[\uu,\vv\right]=\poly\cap L_{-\ww}$ and therefore
$\ww\in L_{\vv}$. Since we also have $-\uu\in L_{\vv}$ and $\ww,-\uu$
are vertices of $\poly$, we get that $\left[\ww,-\uu\right]=\poly\cap L_{\vv}$
is an edge of $\poly$. Therefore, by central symmetry, $\uu,\vv,\ww,-\uu,-\vv,-\ww$
are the vertices of $\poly$.\hfill\halmos
\endproof
\begin{proposition}
\label{prop:six_vertices_area_3}Suppose $\rot\poly=\poly^{\circ}$
and $\left|\verts(\poly)\right|=6$. Then $\vol(\poly)=3$.
\end{proposition}

\proof{Proof.}
Let $\uu,\vv,\ww,-\uu,-\vv,-\ww$ be the six vertices of $\poly$
in the counterclockwise order. We have $\vv,\ww\in L_{\uu}$ and also
$\ww\in L_{\vv}$ by Proposition~\ref{prop:edge_vertex_bijection}.
Thus,
\begin{align*}
\det(\uu,\vv)=\det(\vv,\ww)=\det(\ww,-\uu) & =1
\end{align*}
using the identity $(\rot\xx)^\top\yy=\det(\xx,\yy)$, and hence
\[
\det(-\uu,-\vv)=\det(-\vv,-\ww)=\det(-\ww,\uu)=1.
\]
The sum of these six determinants is equal to twice the area of $\poly$.\hfill\halmos
\endproof

\medskip

\begin{definition}
For a $\zero$-symmetric polygon $\poly$ and a point $\vv\in\R^{2}$, define
\[
\poly_{\vv}:=\conv\left(\poly\cup\left\{ \pm\vv\right\} \right)\cap S_{\vv}.
\]

Note that if $\vv\in\poly^{\circ}$, then $\poly\subseteq S_{\rot\vv}$
which implies
\begin{align}
\rot(\poly_{\rot\vv}) & =\conv\left(\rot\poly\cup\left\{ \pm\vv\right\} \right)\label{eq:rot_P_rot_v-1}\\
(\poly_{\rot\vv})^{\circ} & =\poly^{\circ}\cap S_{\vv}.\label{eq:P_rot_v_polar-1}
\end{align}
\end{definition}

\begin{proposition}[{cf. \cite[Theorem 7.2]{J2021}}]
\label{prop:growing_rot_P}Suppose $\poly$ is a $\zero$-symmetric polygon
such that $\rot\poly\subseteq\poly^{\circ}$. Let $\vv$ be a vertex
of $\poly^{\circ}$. Then
\[
\rot\poly\subseteq\rot(\poly_{\rot\vv})\subseteq(\poly_{\rot\vv})^{\circ}\subseteq\poly^{\circ}.
\]
\end{proposition}

\proof{Proof.}
Equalities (\ref{eq:rot_P_rot_v-1}) and (\ref{eq:P_rot_v_polar-1})
immediately imply the first and third inclusions, respectively. It
remains to show the middle inclusion $\rot(\poly_{\rot\vv})\subseteq(\poly_{\rot\vv})^{\circ}$,
which by (\ref{eq:rot_P_rot_v-1}) and (\ref{eq:P_rot_v_polar-1})
is equivalent to showing
\[
\rot\poly\cup\left\{ \pm\vv\right\} \subseteq\poly^{\circ}\cap S_{\vv}.
\]
We know $\rot\poly\cup\left\{ \pm\vv\right\} \subseteq\poly^{\circ}$
by assumption, so it suffices to show that the strip $S_{\vv}$ contains
both $\pm\vv$ and $\rot\poly$. The strip indeed contains $\pm\vv$
since $(\rot\vv)^\top\vv=0$. To see that the strip contains $\rot\poly$,
observe that $\pm\vv$ are vertices of $\poly^{\circ}$, which means
$\vv^\top\xx=\pm1$ defines two lines spanning parallel edges of $\poly$,
and hence $\pm L_{\vv}$ are two lines spanning parallel edges of
$\rot\poly$. But these two lines form the boundary of the strip $S_{\vv}$.\hfill\halmos
\endproof
\begin{proposition}[{cf. \cite[Corollary 7.3]{J2021}}]
\label{prop:grow_until_equality}Suppose $\poly$ is a symmetric polygon
such that $\rot\poly\subseteq\poly^{\circ}$. Then there exists a
polygon $\qoly$ such that
\[
\rot\poly\subseteq\rot\qoly=\qoly^{\circ}\subseteq\poly^{\circ}.
\]
\end{proposition}

\proof{Proof.}
Assume $\poly^{\circ}\neq\rot\poly$. By Proposition~\ref{prop:growing_rot_P},
it suffices to show that there exists a vertex $\vv$ of $\poly^{\circ}$
not contained in $\rot\poly$ such that
\begin{equation}
\verts(\rot(\poly_{\rot\vv}))\backslash\verts((\poly_{\rot\vv})^{\circ})\subsetneq\verts(\rot\poly)\backslash\verts(\poly^{\circ}).\label{eq:vertex_difference_set_decreases}
\end{equation}
The result then follows by induction on $\left|\verts(\rot\poly)\backslash\verts(\poly^{\circ})\right|$.
After possibly scaling, we may assume without loss of generality that
$\rot\poly$ intersects the boundary of $\poly^{\circ}$.

Assume $\poly^{\circ}\neq\rot\poly$, and let $\vv\in\verts(\poly^{\circ})\backslash\rot\poly$.
We may assume that $\vv$ is contained in an edge of $\poly^{\circ}$
which intersects $\rot\poly$. Indeed, if no such $\vv$ exists, then
we must have that every edge of $\poly^{\circ}$ is contained in $\rot\poly$
or is disjoint from $\rot\poly$. Since $\rot\poly$ intersects the
boundary of $\poly^{\circ}$, this is only possible if the boundaries
of $\rot\poly$ and $\poly^{\circ}$ agree, which contradicts $\poly^{\circ}\neq\rot\poly$.

We start with the containment of (\ref{eq:vertex_difference_set_decreases}).
Let $\ww\in\verts(\rot(\poly_{\rot\vv}))\backslash\verts((\poly_{\rot\vv})^{\circ})$.
Since $(\rot\vv)^\top\vv=0$, it follows by (\ref{eq:P_rot_v_polar-1})
that $\pm\vv$ are vertices of $(\poly_{\rot\vv})^{\circ}$, and hence
we cannot have $\ww=\pm\vv$. Thus $\ww\in V(\rot\poly)$ by (\ref{eq:rot_P_rot_v-1}).
To see that $\ww\notin\verts(\poly^{\circ})$, observe that by Proposition~\ref{prop:growing_rot_P}
we have $\ww\in\rot(\poly_{\rot\vv})\subseteq(\poly_{\rot\vv})^{\circ}\subseteq S_{\vv}$.
However, if it were the case that $\ww\in\verts(\poly^{\circ})$,
then by assumption we would have $\ww\in\verts(\poly^{\circ})\backslash\verts((\poly_{\rot\vv})^{\circ})$,
and by (\ref{eq:P_rot_v_polar-1}) this would imply $\ww\notin S_{\vv}$.

It remains to show the containment of (\ref{eq:vertex_difference_set_decreases})
is strict. Let $\left[\uu,\vv\right]$ be an edge of $\poly^{\circ}$
which intersects $\rot\poly$. If this intersection is given by a
line segment $\left[\rr,\ssv\right]$, so that $\rr$ and $\ssv$ are
distinct vertices of $\rot\poly$ with $\ssv$ a proper convex combination
of $\rr$ and $\vv$, then we have $\ssv\in\verts(\rot\poly)\backslash\verts(\poly^{\circ})$.
On the other hand, by (\ref{eq:rot_P_rot_v-1}) we have $\ssv\notin\verts(\rot(\poly_{\rot\vv}))$
and hence $\ssv\notin\verts(\rot(\poly_{\rot\vv}))\backslash\verts((\poly_{\rot\vv})^{\circ})$.

Otherwise, $\left[\uu,\vv\right]$ intersects $\rot\poly$ at a single
point $\qq\in\verts(\rot\poly)$. In this case, let $\pp,\rr$ be
the neighbouring vertices of $\qq$ in $\rot\poly$. There exists
a unique vertex $\yy$ of $\rot\poly$ such that $\left[\uu,\vv\right]=\poly^{\circ}\cap L_{\yy}$.
Since $\rot\poly\cap L_{\yy}=\left\{ \qq\right\} $, we have that
$\pp,\rr$ lie in the interior of $H_{\yy}$. The edge of $\poly^{\circ}$
spanned by $L_{-\qq}$ contains $\yy$, and the two endpoints $\hat{\uu},\hat{\vv}$
of this edge have $\pp\in L_{\hat{\uu}}$ and $\rr\in L_{\hat{\vv}}$.
Since $L_{\yy}$ is the boundary of $H_{\yy}$, $\yy$ is equal to
neither $\hat{\uu}$ nor $\hat{\vv}$. It follows $\yy\notin\verts(\poly^{\circ})$,
so that $\yy\in\verts(\rot\poly)\backslash\verts(\poly^{\circ})$.

It remains to show $\yy\notin\verts(\rot(\poly_{\rot\vv}))\backslash\verts((\poly_{\rot\vv})^{\circ})$.
For this it suffices to show $\yy\in\verts((\poly_{\rot\vv})^{\circ})$.
Note that $\yy\in\rot\poly\subseteq\poly^{\circ}$ and $\yy\in L_{-\vv}$
by definition of $\yy$. Thus by (\ref{eq:P_rot_v_polar-1}), $\yy$
is contained in $(\poly_{\rot\vv})^{\circ}$. Since $\vv,\ssv$ are
linearly independent by full-dimensionality, we have that $\yy$ is
a vertex of $S_{\vv}\cap S_{\ssv},$which in turn contains $(\poly_{\rot\vv})^{\circ}$.
Therefore we have that $\yy\in\verts((\poly_{\rot\vv})^{\circ})$
as desired.\hfill\halmos
\endproof

\medskip
\begin{definition}
Suppose $\poly$ is a polygon. Let $\stell(\poly)$ denote the set of all points in $\R^{2}$ that
violate at most one inequality constraint of $\poly$. When $\rot\poly=\poly^{\circ}$,
this is
\[
\stell(\poly)=\bigcup_{\vv\in V(\poly)}\bigcap_{\ww\in\verts(\poly)\backslash\{\vv\}}H_{\ww}
\]
\end{definition}

\begin{proposition}
\label{prop:stell_components}Suppose $\rot\poly=\poly^{\circ}$ is
a polygon. Then $\stell(\poly)$ is bounded. In particular, the set
of components of $\stell(\poly)\backslash\poly$ is in bijection with
the set of vertices of $\poly$ as follows:
\begin{align*}
\verts(\poly) & \longrightarrow\left\{ \text{components of }\stell(\poly)\backslash\poly\right\} \\
\vv & \longmapsto H_{\uu}\cap H_{\vv}^{\comp}\cap H_{\ww}
\end{align*}
where, in the above expression, $\uu$ and $\ww$ are the neighbouring
vertices of $\vv$.
\end{proposition}

\proof{Proof.}
By Proposition~\ref{prop:at_least_six_vertices} and central symmetry,
there exists three vertices $\uu,\vv,\ww$ of $\poly$ which are pairwise
linearly independent. We have
\begin{align*}
\stell(\poly) & \subseteq(S_{\uu}\cup S_{\vv})\cap(S_{\uu}\cup S_{\ww})\cap(S_{\vv}\cup S_{\ww})\\
 & =(S_{\uu}\cap S_{\vv})\cup(S_{\uu}\cap S_{\ww})\cup(S_{\vv}\cap S_{\ww})
\end{align*}
which is a union of three parallelograms. So $\stell(\poly)$ is bounded.
The set $\stell(\poly)\backslash\poly$ consists of all points in
the plane which violate exactly one inequality constraint of $\poly$.
Therefore, we have the disjoint union
\[
\stell(\poly)\backslash\poly=\bigcup_{\vv\in\verts(\poly)}H_{\vv}^{\comp}\cap\stell(\poly).
\]
Let $c(\vv)$ denote the closure of $H_{\vv}^{\comp}\cap\stell(\poly)$.
We have $c(\vv$) is a non-empty polygon for each $\vv\in\verts(\poly)$,
hence the above union is a decomposition into the components of $\stell(\poly)\backslash\poly$.
Now fix $\vv\in\verts(\poly)$, and let $\uu,\ww$ be the neighbouring
vertices of $\vv$. We show that $c(\vv)$ is the triangle bounded
by the lines $L_{\uu},L_{\vv},L_{\ww}$. We have $c(\vv)\cap L_{\vv}=\poly\cap L_{\vv}$,
and we denote this edge of $\poly$ by $\left[\hat{\uu},\hat{\ww}\right]$
so that $\hat{\uu}\in L_{\uu}$ and $\hat{\ww}\in L_{\ww}$. The lines
$L_{\pp}$ over all $\pp\in\verts(\poly)$ cut up $L_{\uu}$ and $L_{\ww}$
each into line segments and two half-lines. Let $\left[\hat{\uu},\hat{\vv}'\right]$
be the segment of $L_{\uu}$ that contains $\hat{\uu}$ but $\hat{\vv}'\notin\poly$.
Note that this is indeed a segment and not a half-line, since $\stell(\poly)$
is bounded. Then $\left[\hat{\uu},\hat{\vv}'\right]$ is an edge of
$c(\vv)$. Similarly, let $\left[\hat{\ww},\hat{\vv}''\right]$ be
the segment of $L_{\ww}$ that contains $\hat{\ww}$ but $\hat{\vv}''\notin\poly$.
Here too we have $\left[\hat{\ww},\hat{\vv}''\right]$ is an edge
of $c(\vv)$. 

It remains to show $\hat{\vv}'=\hat{\vv}''$. If these two points
are distinct, then there exists some $\qq$ in $\verts(\poly)$ such
that $\vv,\qq$ are linearly independent and $L_{\qq}$ bounds an
edge of $c(\vv)$ that is disjoint from $\left[\hat{\uu},\hat{\ww}\right]$.
We may assume without loss of generality $\hat{\vv}'\in L_{\qq}$.
Let $\bb\in L_{\qq}\cap L_{\vv}$, and let $\ell$ denote the half-line
\[
\ell=\left\{ \bb+\lambda(\bb-\hat{\vv}'):\lambda\geq0\right\} \subseteq L_{\qq}.
\]
We have $\hat{\vv}'\in H_{\vv}^{\comp}$ while $\bb\in L_{\vv}$,
which implies that $\ell$ contains the edge $\poly\cap L_{\qq}$.
If $\hat{\ww}<\hat{\uu}<\bb$ along $L_{\vv}$, then because $\hat{\ww}\in H_{\uu}\backslash L_{\uu}$
and $\hat{\uu}\in L_{\uu}$, we have $\bb\in H_{\uu}^{\comp}$. Since
$\hat{\vv}'\in L_{\uu}$, $\ell$ does not intersect $\poly$, a contradiction.
Similarly, if $\bb<\hat{\ww}<\hat{\uu}$, then because $\hat{\uu}\in H_{\ww}\backslash L_{\ww}$
and $\hat{\ww}\in L_{\ww}$, we have $\bb\in H_{\ww}^{\comp}$. Since
$\hat{\vv}'\in H_{\ww}$, again we have that $\ell$ does not intersect
$\poly$, a contradiction.\hfill\halmos
\endproof
\begin{proposition}[{cf. \cite[Corollary 7.5]{J2021}}]
\label{prop:P_v_alternative_def}Suppose $\rot\poly=\poly^{\circ}$
is a polygon. Let $\vv\in\stell(\poly)\backslash\poly$. Then
\[
\poly_{\vv}=\conv\left(\left(\poly\cup\left\{ \pm\vv\right\} \right)\cap S_{\vv}\right).
\]
\end{proposition}

\proof{Proof.}
It suffices to show
\[
\conv\left(\poly\cup\left\{ \pm\vv\right\} \right)\backslash\left(\poly\cup\left\{ \pm\vv\right\} \right)\subseteq S_{\vv}.
\]
Choose $\rr\in\conv\left(\poly\cup\left\{ \pm\vv\right\} \right)$
not in $\poly$, not equal to $\pm\vv$. Up to a sign change, $\rr=(1-\lambda)\zz+\lambda\vv$
for some $\zz\in\poly$ and $\lambda\in(0,1)$. There is a unique
vertex $\hat{\vv}\in\verts(\poly)$ such that $L_{\hat{\vv}}$ separates
$\vv$ from $\poly$. We show $\hat{\vv}$ is the unique vertex of
$\poly$ such that $L_{\hat{\vv}}$ separates $\rr$ from $\poly$.
Indeed, if $\rr\in H_{\hat{\rr}}^{\comp}$ for some $\hat{\rr}\in\verts(\poly)$,
then because $\zz\in H_{\hat{\rr}}$ we must have by convexity $\vv\in H_{\hat{\rr}}^{\comp}$.
Therefore, $\hat{\rr}=\hat{\vv}$. Without loss of generality, then,
we assume $\zz\in L_{\hat{\vv}}$. 

Let $L_{\hat{\vv}}\cap\poly=\left[\uu,\ww\right]$ where $\uu,\ww$
are vertices of $\poly$. Note that both $\uu$ and $\ww$ are distinct
from $\hat{\vv}$, one can see this by observing $\uu,\ww\in L_{\hat{\vv}}$
while $(\rot\hat{\vv})^\top\hat{\vv}=0$. Thus $\vv\in S_{\uu}\cap S_{\ww}$,
hence $\uu,\ww\in S_{\vv}$, and so by convexity $\zz\in S_{\vv}$.
Since $\vv\in S_{\vv}$ we get again by convexity $\rr\in S_{\vv}$.\hfill\halmos
\endproof
\begin{proposition}[{cf. \cite[Theorem 7.4]{J2021}}]
\label{prop:P_v_satisfies_rot_equals_polar}Suppose $\rot\poly=\poly^{\circ}$
is a polygon, and $\vv\in\stell(\poly)$. Then we have $\rot\poly_{\vv}=(\poly_{\vv})^{\circ}$.
\end{proposition}

\proof{Proof.}
Since polarity interchanges unions and intersections of closed convex
sets containing the origin, we get by Proposition~\ref{prop:P_v_alternative_def}
that
\begin{align*}
\poly_{\vv}^{\circ} & =\left(\conv\left(\left(\poly\cup\left\{ \pm\vv\right\} \right)\cap S_{\vv}\right)\right)^{\circ}\\
 & =\left(\conv\left(\left(\poly\cap S_{\vv}\right)\cup\left\{ \pm\vv\right\} \right)\right)^{\circ}\\
 & =\left(\poly\cap S_{\vv}\right)^{\circ}\cap(\rot S_{\vv})\\
 & =\conv\left(\poly^{\circ}\cup\left\{ \pm\rot\vv\right\} \right)\cap(\rot S_{\vv})\\
 & =\rot\left(\conv\left(\rot\poly^{\circ}\cup\left\{ \pm\vv\right\} \right)\cap S_{\vv}\right)\\
 & =\rot\left(\conv\left(\poly\cup\left\{ \pm\vv\right\} \right)\cap S_{\vv}\right)\\
 & =\rot\poly_{\vv}.\tag*{\halmos}
\end{align*}
\endproof

The next definition gives names to the points along the boundary of $\poly$ which are involved in the transformation from $\poly$ into $\poly_{\vv}$.
\begin{definition}
\label{def:vertices_along_boundary_of_P}Let $\poly$ be a polygon
such that $\rot\poly=\poly^{\circ}$, and let $\vv\in\stell(\poly)\backslash\poly$.
With dependence on the pair $\left(\poly,\vv\right)$, we define
\[
-\ww\leq\pp\leq\hat{\uu}<\hat{\vv}<\hat{\ww}\leq\qq\leq\uu<\ww
\]
to be the points along the boundary of $\poly$ in the counterclockwise
order such that $\hat{\vv}$ is the unique vertex for which $L_{\hat{\vv}}$
separates $\vv$ from $\poly$, and
\begin{align*}
\poly\cap L_{\hat{\vv}} & =\left[\uu,\ww\right] & \poly\cap L_{-\uu} & =\left[\pp,\hat{\vv}\right]\\
\poly\cap L_{-\vv} & =\left[\hat{\uu},\hat{\ww}\right] & \poly\cap L_{-\ww} & =\left[\hat{\vv},\qq\right].
\end{align*}
\end{definition}

\begin{figure}[t]
\begin{centering}
\includegraphics[scale=0.4]{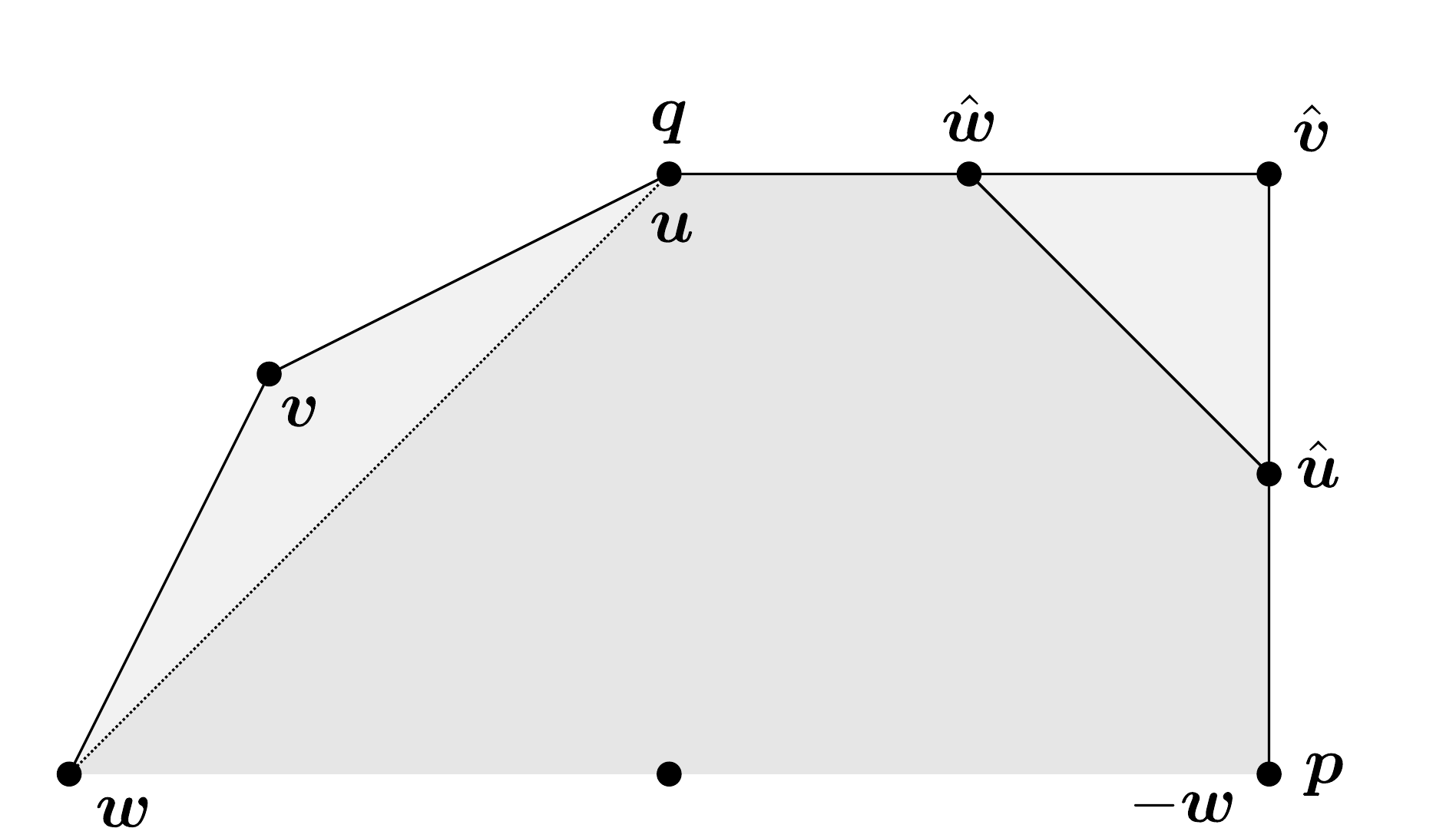}
\par\end{centering}
\caption{\label{fig:vertices_along_P}The boundary points of $\protect\poly$
given in Definition~\ref{def:vertices_along_boundary_of_P}. In this example, $\qq=\uu$ and $\pp=-\ww$, but these equalities need not hold in general.}
\end{figure}

See Figure~\ref{fig:vertices_along_P} for an illustration. Note
that $\pp,\hat{\vv},\qq,\uu,\ww$ are all vertices of $\poly$, and
that $\pp\leq\hat{\uu}<\hat{\vv}<\hat{\ww}\leq\qq$ since $-\vv$
is the unique vertex of $\poly$ for which $L_{-\vv}$ separates $\hat{\vv}$
from all other vertices of $\poly$. We adopt the notation of Definition~\ref{def:vertices_along_boundary_of_P}
in Propositions~\ref{prop:vertex_symm_diff_1} and~\ref{prop:vertex_symm_diff_2}
below.
\begin{proposition}
\label{prop:vertex_symm_diff_1}The symmetric difference of $\verts(\poly_{\vv})$
and $\verts(\poly)$ satisfies
\[
\left\{ \pm\vv,\pm\hat{\vv}\right\} \subseteq\verts(\poly_{\vv})\triangle\verts(\poly)\subseteq\left\{ \pm\vv,\pm\hat{\vv},\pm\uu,\pm\hat{\uu},\pm\ww,\pm\hat{\ww}\right\} .
\]
\end{proposition}

\proof{Proof.}
The boundary of $S_{\vv}$ is $\pm L_{\vv}$, so any vertex of $\poly_{\vv}$
that is not contained in $\poly$ must lie in 
\[
(\poly\cup\left\{ \pm\vv\right\} )\cap(\pm L_{\vv})=\poly\cap(\pm L_{\vv})=\pm\left[\hat{\uu},\hat{\ww}\right]
\]
by Proposition~\ref{prop:P_v_alternative_def}. Hence
\[
\left\{ \pm\vv\right\} \subseteq\verts(\poly_{\vv})\backslash\verts(\poly)\subseteq\left\{ \pm\vv,\pm\hat{\uu},\pm\hat{\ww}\right\} .
\]
We show
\[
\left\{ \pm\hat{\vv}\right\} \subseteq\verts(\poly)\backslash\verts(\poly_{\vv})\subseteq\left\{ \pm\hat{\vv},\pm\uu,\pm\ww\right\} .
\]
By assumption, $\vv\in H_{\hat{\vv}}^{\comp}$ and therefore $-\hat{\vv}\in H_{\vv}^{\comp}$
which shows $\hat{\vv}\in\verts(\poly)\backslash\verts(\poly_{\vv})$.
Now, let $\bb\in\verts(\poly)\backslash\verts(\poly_{\vv})$ distinct
from $\pm\hat{\vv}$. Since $\pm\hat{\uu},\pm\hat{\ww}\in\verts(\poly_{\vv})$,
$\bb$ must be a vertex of $\conv(\poly\cup\left\{ \pm\vv\right\} )$
by Proposition~\ref{prop:P_v_alternative_def}. Let $\rr,\ssv$ be
the vertices of $\poly$ such that $\bb\in L_{\rr}\cap L_{\ssv}$.
If $\vv\in S_{\rr}\cap S_{\ssv}$ then $\bb$ would be a vertex of
$S_{\rr}\cap S_{\ssv}$ which is a parallelogram containing $\conv(\poly\cup\left\{ \pm\vv\right\} )$,
and hence $\bb$ would be a vertex of $\conv(\poly\cup\left\{ \pm\vv\right\} )$,
a contradiction. So either $\vv\notin S_{\rr}$ or $\vv\notin S_{\ssv}$.
It follows that either $\rr=\pm\hat{\vv}$ or $\ssv=\pm\hat{\vv}$.
Therefore, $\bb\in\poly\cap(\pm L_{\hat{\vv}})=\pm\left[\uu,\ww\right]$,
so $\bb\in\left\{ \pm\uu,\pm\ww\right\} $.\hfill\halmos
\endproof
\begin{proposition}
\label{prop:vertex_symm_diff_2}Suppose it is not the case that both
$\left|\verts(\poly)\right|=6$ and $\vv\in L_{\uu}\cap L_{-\ww}$.
Then the symmetric difference of $\verts(\poly_{\vv})$ and $\verts(\poly)$
has size 8 and is given by
\[
\verts(\poly_{\vv})\triangle\verts(\poly)=\left\{ \pm\vv,\pm\hat{\vv},\pm\tilde{\uu},\pm\tilde{\ww}\right\} 
\]
for some $\tilde{\uu}\in\left\{ \pm\uu,\pm\hat{\uu}\right\} $ and
some $\tilde{\ww}\in\left\{ \pm\ww,\pm\hat{\ww}\right\} $.
\end{proposition}

\proof{Proof.}
We have $\uu,\ww\in\verts(\poly)$ and $\hat{\uu},\hat{\ww}\in\verts(\poly_{\vv})$.
Applying Proposition~\ref{prop:vertex_symm_diff_1}, we would like
to show that $\hat{\uu}\in\verts(\poly)$ if and only if $\uu\notin\verts(\poly_{\vv})$,
and similarly $\hat{\ww}\in\verts(\poly)$ if and only if $\ww\notin\verts(\poly_{\vv})$.
We sketch the argument of the former claim; the latter claim is analogous.

Suppose $\hat{\uu}\in\verts(\poly)$. Then $\hat{\uu}=\pp$. This
implies $\uu,\vv\in L_{\pp}$. Suppose for a contradiction $\uu\in\verts(\poly_{\vv})$.
Then $\poly_{\vv}\cap L_{\pp}=\left[\uu,\vv\right]$, and since $\uu\in\poly\cap L_{\hat{\vv}}$
but $L_{\hat{\vv}}$ separates $\vv$ from $\poly$, we get $\poly\cap\left[\uu,\vv\right]=\left\{ \uu\right\} $.
Meanwhile $\poly\cap L_{\pp}=\left[\uu',\uu\right]$ for some vertex
$\uu'$ of $\poly$, which shows $\uu\in(\uu',\vv)$. It follows that
$\uu'$ is a vertex of $\poly$ outside of $S_{\vv}$, since otherwise
$\uu$ would not be a vertex of $\poly_{\vv}$. Since $-\ww<\hat{\vv}\leq\uu'<\uu<\ww$
we have $\uu'=\hat{\vv}$. Hence $\uu=\qq$. It follows that $\hat{\vv}<\uu<\ww$
are consecutive vertices of $\poly$ in the counterclockwise ordering.
Since $\uu,\ww\in L_{\hat{\vv}}$, we get by Proposition~\ref{prop:sufficient_condition_for_six_verts}
that $\left|\verts(\poly)\right|=6$. We have $\vv\in L_{-\ww}$ since
there are only six vertices and hence $-\ww=\pp=\hat{\uu}$. Since
$\poly_{\vv}\cap L_{\pp}=\left[\uu,\vv\right]$ and since $\uu\in(\uu',\vv)$,
we have $\poly_{\vv}\cap\left[\uu',\uu\right]=\left\{ \uu\right\} $.
Since $\hat{\ww}$ is a vertex of $\poly_{\vv}$ for which $\hat{\ww}\in\left[\hat{\vv},\qq\right]=\left[\uu',\uu\right]$,
we get $\uu=\hat{\ww}$ and therefore $\vv\in L_{\uu}$. Thus $\vv\in L_{\uu}\cap L_{-\ww}$.
This contradicts the given hypotheses of the proposition. We conclude
$\uu\notin\verts(\poly_{\vv})$.

Now suppose $\hat{\uu}\notin\verts(\poly)$. Since $L_{\hat{\vv}}$
separates $\vv$ from all other vertices of $\poly$, we have that
$L_{-\vv}$ separates $\hat{\vv}$ from all other vertices of $\poly$.
Therefore, since $\poly\cap L_{-\vv}=\left[\hat{\uu},\hat{\ww}\right]$,
we have $\hat{\uu}$ and $\hat{\ww}$ both lie in edges incident to
$\hat{\vv}$. Since $\pp\leq\hat{\uu}<\hat{\vv}$ we have in particular
that $\hat{\uu}$ lies in the edge $\left[\pp,\hat{\vv}\right]$ of
$\poly$. But this edge is given by the intersection $\poly\cap L_{-\uu}$,
which implies $\hat{\uu}\in L_{-\uu}$, and therefore $\uu\in L_{\hat{\uu}}$.
We also have $\uu\in L_{\pp}$, and also $\uu\in\poly_{\vv}$ since
$\uu\in S_{\vv}$. Now $\pp\in\verts(\poly)$, and by assumption,
$\hat{\uu}\notin\verts(\poly)$, which implies $\hat{\uu}\neq\pp$.
Both $\hat{\uu},\pp$ are vertices of $\poly_{\vv}$, the latter since
$\pp\neq\hat{\vv}$. Therefore, since $\uu\in\poly_{\vv}$ is contained
in the intersection of the two edges $\poly_{\vv}\cap L_{\pp}$ and
$\poly_{\vv}\cap L_{\hat{\uu}}$, we conclude $\uu\in\verts(\poly_{\vv})$.\hfill\halmos
\endproof
\begin{proposition}
\label{prop:Q_equals_P_v}Suppose $\rot\qoly=\qoly^{\circ}$ is a
polygon such that $\left|\verts(\qoly)\right|>6$, and let $\vv\in\verts(\qoly)$.
Then there exists a polygon $\poly$ such that $\rot\poly=\poly^{\circ}$,
$\vv$ lies in the interior of $\stell(\poly)\backslash\poly$, and
$\qoly=\poly_{\vv}$.
\end{proposition}

\proof{Proof.}
Let $\uu,\ww$ be the vertices of $\qoly$ adjacent to $\vv$, so
that $\uu<\vv<\ww$ in the counterclockwise order. The triangle bounded
by $L_{-\uu},L_{-\vv},L_{-\ww}$ is the closure of a component of
$\stell(\qoly)\backslash\qoly$ by Proposition~\ref{prop:stell_components}.
Let $\hat{\vv}$ be the unique point of $L_{-\uu}\cap L_{-\ww}$,
which is the unique vertex of this triangle not in $\qoly$. Let $\hat{\uu}<\hat{\ww}\in\verts(\qoly)$
be the other two vertices of this triangle, so that $\qoly\cap L_{-\vv}=\left[\hat{\uu},\hat{\ww}\right]$,
and let $\poly=\qoly_{\hat{\vv}}$. We have
\begin{equation}
\left[\uu,\ww\right]=\qoly\cap L_{\hat{\vv}}=\poly\cap L_{\hat{\vv}}\label{eq:P_cap_L_hat_v}
\end{equation}
where the second equality holds by Proposition~\ref{prop:P_v_alternative_def}.
Since $\rot\poly=\poly^{\circ}$ by Proposition~\ref{prop:P_v_satisfies_rot_equals_polar},
we have $\uu,\ww$ are vertices of both $\poly$ and $\qoly$. Since
$\left|\verts(\qoly)\right|>6$, we therefore apply Proposition~\ref{prop:vertex_symm_diff_2}
with respect to the pair $\left(\qoly,\hat{\vv}\right)$, which corresponds
to the Definition~\ref{def:vertices_along_boundary_of_P} sequence
\[
-\hat{\ww}\leq-\uu\leq-\uu<-\vv<-\ww\leq-\ww\leq\hat{\uu}<\hat{\ww}
\]
of boundary points of $\qoly$, to get
\begin{equation}
\verts(\poly)\triangle\verts(\qoly)=\left\{ \pm\vv,\pm\hat{\vv},\pm\hat{\uu},\pm\hat{\ww}\right\} .\label{eq:vertices_P_Q_symm_diff}
\end{equation}

Let $\pp<\hat{\vv}<\qq$ be consecutive vertices of $\poly$ in the
counterclockwise ordering. Then $L_{\pp},L_{\hat{\vv}},L_{\qq}$ bound
the closure of a component of $\stell(\poly)\backslash\poly$ by Proposition~\ref{prop:stell_components}.
We show $\vv$ lies in the interior of this component. The fact that
$\hat{\vv}\notin H_{-\vv}$ implies $\vv\notin H_{\hat{\vv}}$. To
see $\vv\in H_{\pp}\cap H_{\qq}$ it suffices to show $\pp,\qq$ are
vertices of $\qoly$, which is equivalent to saying $\pp,\qq\notin\verts(\poly)\triangle\verts(\qoly)$.
Since $\vv\notin\poly$ and $\pp<\hat{\vv}<\qq$ we have $\pp,\qq\notin\left\{ \pm\vv,\pm\hat{\vv}\right\} $.
Since $\hat{\uu},\hat{\ww}$ are vertices of $\qoly$, by (\ref{eq:vertices_P_Q_symm_diff})
they are not vertices of $\poly$, and hence $\pp,\qq\notin\left\{ \pm\hat{\uu},\pm\hat{\ww}\right\} $.
We also have $\vv\notin L_{\pp}$ and $\vv\notin L_{\qq}$, since
otherwise we would have $\vv\in\qoly\cap L_{\hat{\uu}}\cap L_{\hat{\ww}}\cap L_{\pp}$
or $\vv\in\qoly\cap L_{\hat{\uu}}\cap L_{\hat{\ww}}\cap L_{\qq}$.
As $\hat{\uu},\hat{\ww},\pp,\qq\in\verts(\poly)$ are pairwise distinct,
either case would imply that $\vv$ lies in the intersection of three
edges of $\qoly$, a contradiction.

The final step is to show $\poly_{\vv}=\qoly$. We begin by showing
$\poly\cap L_{\hat{\vv}}=\left[\uu,\ww\right]$ and $\poly\cap L_{-\vv}=\left[\hat{\uu},\hat{\ww}\right]$.
The former equality holds by (\ref{eq:P_cap_L_hat_v}). We establish
the latter equality. Observe that $\hat{\uu},\hat{\ww}$ both lie
on the boundary of $\poly$. Indeed, $\hat{\uu},\hat{\ww}$ are vertices
of $\qoly$ distinct from $\vv$, which means $\hat{\uu},\hat{\ww}\in S_{\hat{\vv}}$
and therefore $\hat{\uu},\hat{\ww}\in\poly$. Since $\uu,\ww\in\verts(\poly)$,
we have $\poly\cap L_{-\uu}$ and $\poly\cap L_{-\ww}$ are two edges
of $\poly$ which contain $\hat{\uu}$ and $\hat{\ww}$, respectively.
The fact that $\hat{\uu},\hat{\ww}\in L_{-\vv}$ concludes the claim
$\poly\cap L_{-\vv}=\left[\hat{\uu},\hat{\ww}\right]$.

Note that $\vv\notin L_{\uu}\cap L_{-\ww}$. This is equivalent to
the statement $\qoly\cap L_{-\vv}\neq\left[\uu,-\ww\right]$, which
is true because we already know $\qoly\cap L_{-\vv}=\left[\hat{\uu},\hat{\ww}\right]$
and that $\hat{\uu},\hat{\ww}\notin\verts(\poly)$ by (\ref{eq:vertices_P_Q_symm_diff}).
We have
\[
\left[\hat{\uu},\hat{\ww}\right]=\poly\cap L_{-\vv}=\poly_{\vv}\cap L_{-\vv}
\]
where again the second equality holds by Proposition~\ref{prop:P_v_alternative_def},
and since $-\vv\in\verts(\poly_{\vv})$ and $\rot\poly_{\vv}=(\poly_{\vv})^{\circ}$
we get as before $\hat{\uu},\hat{\ww}\in\verts(\poly_{\vv})$. We
again apply Proposition~\ref{prop:vertex_symm_diff_2}, this time
in terms of the pair $\left(\poly,\vv\right)$. We get
\[
\verts(\poly_{\vv})\triangle\verts(\poly)=\left\{ \pm\vv,\pm\hat{\vv},\pm\hat{\uu},\pm\hat{\ww}\right\} .
\]
Thus the vertex sets of $\poly_{\vv}$ and $\qoly$ agree.\hfill\halmos
\endproof

\medskip

\begin{definition}
\label{def:interpolated_P_v}Let $\qoly=\poly_{\vv}$ as in Proposition~\ref{prop:Q_equals_P_v},
and assume the notation of Definition~\ref{def:vertices_along_boundary_of_P}
with respect to the pair $(\poly,\vv)$. Let $\vv_{0}=\uu$ and let
$\vv_{1}\in L_{\hat{\uu}}\cap L_{\qq}$. For $\lambda\in\left[0,1\right]$,
let $\vv_{\lambda}=(1-\lambda)\vv_{0}+\lambda\vv_{1}$. Let $\hat{\ww}_{\lambda}$
denote the unique point in $L_{-\vv_{\lambda}}\cap L_{-\ww}$. Finally,
let
\[
\poly_{\lambda}:=\poly_{\vv_{\lambda}}.
\]
\end{definition}

Note that this definition is well-defined since $\vv_{\lambda}\in\stell(\poly)$.
By Proposition~\ref{prop:P_v_satisfies_rot_equals_polar}, we have
$\rot\poly_{\lambda}=(\poly_{\lambda})^{\circ}$. Note that $\poly=\poly_{0}$
and that $\vv=\vv_{\lambda}$ for some $\lambda\in(0,1)$ since $\vv\in L_{\hat{\uu}}\cap L_{\hat{\ww}}$ and $L_{\hat{\uu}}$ intersects the boundary of the component of $\stell(\poly)\backslash\poly$ containing $\vv$ at $\vv_0$ and $\vv_1$. Hence, we have $\qoly=\poly_{\vv}=\poly_{\vv_{\lambda}}=\poly_{\lambda}$. See Figure~\ref{fig:interpolated_P_v} for an illustration of Definition~\ref{def:interpolated_P_v}.

\begin{figure}[t]
\begin{centering}
\includegraphics[scale=0.4]{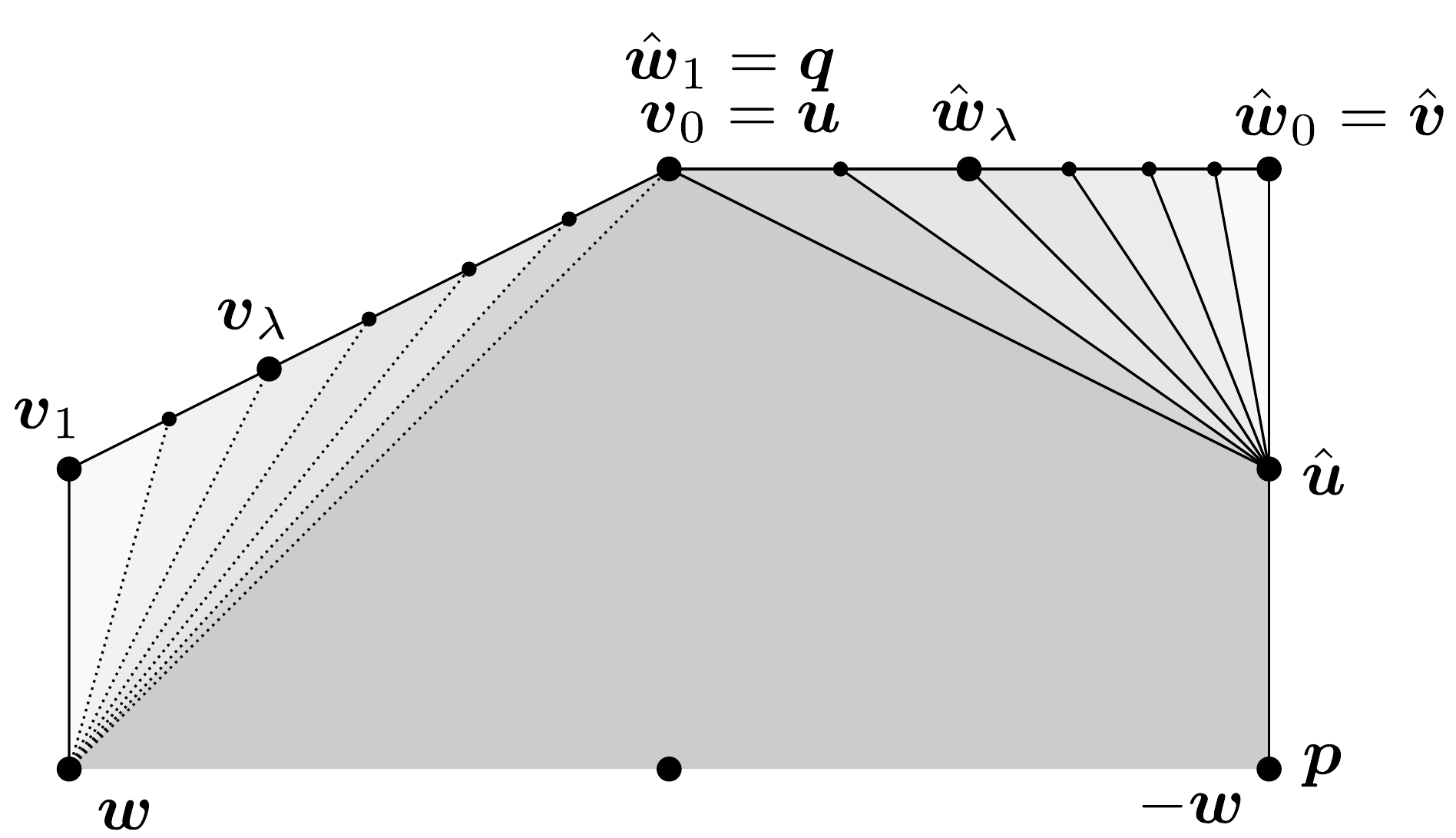}
\par\end{centering}
\caption{\label{fig:interpolated_P_v}The polytopes $\protect\poly_{\lambda}$
for $\lambda\in\left\{ 0,\tfrac{1}{6},\frac{2}{6},\tfrac{3}{6},\frac{4}{6},\frac{5}{6},1\right\} $.}
\end{figure}
\begin{proposition}
\label{prop:P_lambda_mu_coefficient}Let $\lambda\in\left[0,1\right]$,
and write $\hat{\ww}_{\lambda}=(1-\mu)\hat{\vv}+\mu\qq$. Then
\[
\mu=\frac{\lambda a}{\lambda a+(1-\lambda)b},
\]
where $a:=\det(\hat{\vv},\vv_{1}-\uu)$ and $b:=\det(\uu,\qq-\hat{\vv})$.
\end{proposition}
Figure~\ref{fig:interpolated_P_v} demonstrates the nonlinear dependence of $\mu$ on $\lambda$.
\proof{Proof.}
Since $\hat{\ww}_{\lambda}\in L_{-\vv_{\lambda}}$ we have
\begin{align*}
1 & =-(\rot\vv_{\lambda})^\top\hat{\ww}_{\lambda}\\
 & =\det((1-\mu)\hat{\vv}+\mu\qq,(1-\lambda)\uu+\lambda\vv_{1})\\
 & =\det(\hat{\vv}+\mu(\qq-\hat{\vv}),\uu+\lambda(\vv_{1}-\uu))\\
 & =1+\lambda a-\mu b+\lambda\mu\det(\qq-\hat{\vv},\vv_{1}-\uu).
\end{align*}
Since $\uu\in L_{\hat{\vv}}$ and $\vv_{1}\in L_{\qq}$ we have $\det(\qq,\vv_{1})=1=\det(\hat{\vv},\uu)$.
Hence
\begin{align}
\det(\qq-\hat{\vv},\vv_{1}-\uu) & =\det(\qq,\vv_{1})-\det(\qq,\uu)-a\label{eq:beta_minus_alpha}\\
 & =\det(\hat{\vv},\uu)-\det(\qq,\uu)-a\nonumber \\
 & =b-a,\nonumber 
\end{align}
and therefore we get
\[
1=1+\lambda a-\mu b+\lambda\mu(b-a).
\]
Solving for $\mu$ yields the desired equality.\hfill\halmos
\endproof
\begin{proposition}
\label{prop:P_lambda_area_min_attained_on_boundary}For $\lambda\in\left[0,1\right]$,
we have
\begin{align*}
\vol(\poly_{\lambda}) & \geq\min\{\vol(\poly_{0}),\vol(\poly_{1})\}.
\end{align*}
\end{proposition}

\proof{Proof.}
Observe that
\begin{align*}
\vol(\poly_{0}\backslash\poly_{\lambda}) & =\left|\det(\hat{\vv}-\hat{\uu},\hat{\ww}_{\lambda}-\hat{\uu})\right|=\mu\left|\det(\hat{\vv}-\hat{\uu},\qq-\hat{\uu})\right|\\
\vol(\poly_{\lambda}\backslash\poly_{0}) & =\left|\det(\uu-\ww,\vv_{\lambda}-\ww)\right|=\lambda\left|\det(\uu-\ww,\vv_{1}-\ww)\right|.
\end{align*}
We have
\[
\vol(\poly_{\lambda})=\vol(\poly_{0})+\vol(\poly_{\lambda}\backslash\poly_{0})-\vol(\poly_{0}\backslash\poly_{\lambda}),
\]
and so
\[
\frac{\mathrm{d}^{2}\vol(\poly_{\lambda})}{\mathrm{d}\lambda^{2}}=-\left|\det(\hat{\vv}-\hat{\uu},\qq-\hat{\uu})\right|\cdot\frac{\mathrm{d}^{2}\mu}{\mathrm{d}\lambda^{2}}.
\]
Hence we are done if we can show $\mu=\mu(\lambda)$ is convex on
$\lambda\in\left[0,1\right]$, as this would imply that the minimum
of $\vol(\poly_{\lambda})$ is attained at either $\lambda=0$ or
$\lambda=1$. By Proposition~\ref{prop:P_lambda_mu_coefficient},
we have
\[
\frac{\mathrm{d}^{2}\mu}{\mathrm{d}\lambda^{2}}=\frac{2ab(b-a)}{(\lambda a+(1-\lambda)b)^{3}}.
\]
It therefore remains to show $b\geq a>0$. Since $\vv_{1}$ is separated
from $\poly$ by $L_{\hat{\vv}}$ we have $\det(\hat{\vv},\vv_{1})>1$.
Since $\det(\hat{\vv},\uu)=1$ we get $a>0$. To see that $b-a>0$,
we use the representation of (\ref{eq:beta_minus_alpha}) to write
$b-a=\det(\qq-\hat{\vv},\vv_{1}-\uu)$. We have $\qq,\hat{\vv}\in L_{-\ww}$,
which implies $\qq-\hat{\vv}$ is a scalar multiple of $\ww$. Since
$-\ww<\hat{\vv}<\qq<\ww$ along the boundary of $\poly$ in the counterclockwise
order, we have that $\qq-\hat{\vv}$ is a positive multiple of $\ww$.
In a similar manner, we have $\vv_{1},\uu\in L_{\hat{\uu}}$ which
implies $\vv_{1}-\uu$ is a scalar multiple of $\hat{\uu}$. Since
$-\ww<\hat{\uu}<\uu<\vv_{1}<\ww<-\hat{\uu}$ along the boundary of
$\poly_{1}$ in the counterclockwise order, $\vv_{1}-\uu$ is a negative
multiple of $\hat{\uu}$. We conclude that $b-a$ has the same sign
as $\det(\ww,-\hat{\uu})$. Since $-\ww<\hat{\uu}<\ww$ in the counterclockwise
order of $\poly$, this determinant is positive.\hfill\halmos
\endproof
\begin{proposition}
\label{prop:P_lambda_fewer_verts_on_boundary}For $\lambda\in(0,1)$,
we have
\begin{align*}
\left|\verts(\poly_{\lambda})\right| & >\max\{\left|\verts(\poly_{0})\right|,\left|\verts(\poly_{1})\right|\}.
\end{align*}
\end{proposition}

\proof{Proof.}
We have $\hat{\uu}\in(\pp,\hat{\vv})$ and therefore $\hat{\uu}\notin\verts(\poly)$.
Since $\lambda\in(0,1)$, and $\hat{\ww}_{\lambda}\in(\hat{\vv},\qq)$,
we also have $\hat{\ww}_{\lambda}\notin\verts(\poly)$. Therefore,
by Proposition~\ref{prop:vertex_symm_diff_2}, we have
\[
\verts(\poly_{\lambda})=(\verts(\poly)\backslash\left\{ \pm\hat{\vv}\right\} )\cup\left\{ \pm\vv_{\lambda},\pm\hat{\uu},\pm\hat{\ww}_{\lambda}\right\} .
\]
Now $L_{-\vv_{1}}\cap\poly=\left[\hat{\uu},\qq\right]$. We have $\qq\in\verts(\poly)$,
and therefore by Proposition~\ref{prop:vertex_symm_diff_2}, we have
\[
\verts(\poly_{1})=(\verts(\poly)\backslash\left\{ \pm\hat{\vv},\pm\qq\right\} )\cup\left\{ \pm\vv_{1},\pm\hat{\uu}\right\} .
\]
Since $\vv_{0}=\uu$, we have $\poly_{0}=\poly$, and therefore we
conclude
\[
\left|\verts(\poly_{\lambda})\right|>\left|\verts(\poly_{0})\right|=\left|\verts(\poly_{1})\right|.\tag*{\halmos}
\]
\endproof
Recall the statement of Lemma~\ref{lem:polar_area_lower_bound}:
if $\qoly$ is a polygon satisfying $\rot\qoly\subseteq\qoly^{\circ}$
then $\vol(\qoly^{\circ})\geq3$.

\medskip
\proof{Proof of Lemma~\ref{lem:polar_area_lower_bound}.}
Suppose $\qoly$ is a polygon satisfying $\rot\qoly\subseteq\qoly^{\circ}$.
By Proposition~\ref{prop:grow_until_equality}, we may assume without
loss of generality that $\rot\qoly=\qoly^{\circ}$. Then $\left|\verts(\qoly)\right|\geq6$
by Proposition~\ref{prop:at_least_six_vertices}. If $\left|\verts(\qoly)\right|=6$
then $\vol(\qoly)=3$ by Proposition~\ref{prop:six_vertices_area_3}.
Otherwise, $\left|\verts(\qoly)\right|>6$. By Proposition~\ref{prop:Q_equals_P_v},
there exists $\vv\in\verts(\qoly)$ such that $\qoly=\poly_{\vv}$
for some $\vv$ in the interior of $\stell(\poly)\backslash\poly$.
For $\lambda\in\left[0,1\right]$, let $\poly_{\lambda}$ be the polytope
of Definition~\ref{def:interpolated_P_v}, in terms of the pair $(\poly,\vv)$,
so that in particular there exists some $\lambda\in(0,1)$ such that
$\qoly=\poly_{\lambda}$. By Proposition~\ref{prop:P_lambda_area_min_attained_on_boundary},
there exists $i\in\left\{ 0,1\right\} $ for which $\vol(\qoly)\geq\vol(\poly_{i})$.
By Proposition~\ref{prop:P_lambda_fewer_verts_on_boundary}, $\left|\verts(\qoly)\right|>\left|\verts(\poly_{i})\right|$.
By induction on the number of vertices, we have $\vol(\poly_{i})\geq3$,
and therefore $\vol(\qoly^{\circ})=\vol(\rot\qoly)=\vol(\qoly)\geq3$.
\hfill\halmos
\endproof

\end{APPENDICES}

\subsubsection*{Acknowledgements.} J. Paat was supported by a Natural Sciences and Engineering Research Council of Canada (NSERC) Discovery Grant [RGPIN-2021-02475]. 
R. Weismantel was supported by the Einstein Foundation Berlin. 

\bibliographystyle{informs2014}
\bibliography{manuscript}

\begin{thebibliography}{32}
\providecommand{\natexlab}[1]{#1}
\providecommand{\url}[1]{\texttt{#1}}
\providecommand{\urlprefix}{URL }

\bibitem[{Aliev et~al.(2021)Aliev, Celaya, Henk, \protect\BIBand{}
  Williams}]{ACHW2021}
Aliev I, Celaya M, Henk M, Williams A (2021) Distance-sparsity transference for
  vertices of corner polyhedra. \emph{{SIAM} Journal on Optimization}
  31:200--126, \urlprefix\url{http://dx.doi.org/10.1137/20M1353228}.

\bibitem[{Aliev et~al.(2020)Aliev, Henk, \protect\BIBand{} Oertel}]{AHO2019}
Aliev I, Henk M, Oertel T (2020) Distances to lattice points in knapsack
  polyhedra. \emph{Mathematical Programming} 182:175--198,
  \urlprefix\url{http://dx.doi.org/10.1007/s10107-019-01392-1}.

\bibitem[{Celaya et~al.(2022)Celaya, Kuhlmann, Paat, \protect\BIBand{}
  Weismantel}]{CKPW2022}
Celaya M, Kuhlmann S, Paat J, Weismantel R (2022) Improving the {C}ook et al.
  proximity bound given integral valued constraints. \emph{Integer Programming
  and Combinatorial Optimization}, 84 -- 97 (Springer International
  Publishing), \urlprefix\url{http://dx.doi.org/10.1007/978-3-031-06901-7\_7}.

\bibitem[{Cook et~al.(1986)Cook, Gerards, Schrijver, \protect\BIBand{}
  Tardos}]{CGST1986}
Cook W, Gerards A, Schrijver A, Tardos E (1986) Sensitivity theorems in integer
  linear programming. \emph{Mathematical Programming} 34:251--264,
  \urlprefix\url{http://dx.doi.org/10.1007/BF01582230}.

\bibitem[{{Del Pia} \protect\BIBand{} Ma(2021)}]{DM2021}
{Del Pia} A, Ma M (2021) Proximity in concave integer quadratic programming.
  \emph{Available online: arXiv:2006.01718} .

\bibitem[{Eisenbrand \protect\BIBand{} Weismantel(2020)}]{EW2018}
Eisenbrand F, Weismantel R (2020) Proximity {R}esults and {F}aster {A}lgorithms
  for {I}nteger {P}rogramming {U}sing the {S}teinitz {L}emma. \emph{{ACM}
  Transactions on Algorithms} 16:1--14,
  \urlprefix\url{http://dx.doi.org/10.1145/3340322}.

\bibitem[{Fischetti et~al.(2005)Fischetti, Glover, \protect\BIBand{}
  Lodi}]{FGL2005}
Fischetti M, Glover F, Lodi A (2005) The feasibility pump. \emph{Mathematical
  Programming} 104:91--104,
  \urlprefix\url{http://dx.doi.org/10.1007/s10107-004-0570-3}.

\bibitem[{Florian(1996)}]{F1996}
Florian A (1996) On the area sum of a convex set and its polar reciprocal.
  \emph{Mathematica Pannonica} 171:176.

\bibitem[{Fortier(2020)}]{F2020}
Fortier JM (2020) \emph{Self-nolar Planar Polytopes: When Finding the Polar is
  Rotating by Pi}. Master's thesis, Concordia University.

\bibitem[{Granot \protect\BIBand{} Skorin-Kapov(1990)}]{G1990}
Granot F, Skorin-Kapov J (1990) Some proximity and sensitivity results in
  quadratic integer programming. \emph{Mathematical Programming} 47:259--268,
  \urlprefix\url{http://dx.doi.org/10.1007/BF01580862}.

\bibitem[{Gribanov \protect\BIBand{} Veselov(2016)}]{GV2016}
Gribanov D, Veselov S (2016) On integer programming with bounded determinants.
  \emph{Optimization Letters} 10:1169--1177.

\bibitem[{Gruber(2007)}]{G2007}
Gruber P (2007) \emph{Convex and Discrete Geometry} (Springer-{V}erlag {B}erlin
  {H}eidelberg).

\bibitem[{Henk et~al.(2022)Henk, Kuhlmann, \protect\BIBand{}
  Weismantel}]{HKW2022}
Henk M, Kuhlmann S, Weismantel R (2022) On lattice width of lattice-free
  polyhedra and height of {H}ilbert bases. \emph{SIAM Journal on Discrete
  Mathematics} 36(3):1918–1942.

\bibitem[{Hochbaum \protect\BIBand{} Shanthikumar(1990)}]{H1990}
Hochbaum DS, Shanthikumar JG (1990) Convex separable optimization is not much
  harder than linear optimization. \emph{Journal of the ACM} 37:843--862,
  \urlprefix\url{http://dx.doi.org/10.1145/96559.96597}.

\bibitem[{Jansen \protect\BIBand{} Rohwedder(2019)}]{JR2018}
Jansen K, Rohwedder L (2019) {On Integer Programming and Convolution}.
  \emph{10th Innovations in Theoretical Computer Science}, volume~43,
  43:1--43:17, \urlprefix\url{http://dx.doi.org/10.4230/LIPIcs.ITCS.2019.43}.

\bibitem[{Jensen(2021)}]{J2021}
Jensen A (2021) Self-polar polytopes. \emph{Polytopes and Discrete Geometry},
  Contemporary Mathematics (American Mathematical Society),
  \urlprefix\url{http://dx.doi.org/10.1090/conm/764}.

\bibitem[{Jiang \protect\BIBand{} Basu(2022)}]{BJ2022}
Jiang H, Basu A (2022) Enumerating integer points in polytopes with bounded
  subdeterminants. \emph{SIAM Journal on Discrete Mathematics} 36(1):449--460,
  \urlprefix\url{http://dx.doi.org/10.1137/21M139935X}.

\bibitem[{Khinchine(1948)}]{K1948}
Khinchine A (1948) A quantitative formulation of {K}ronecker's theory of
  approximation. \emph{Izv. Acad. Nauk SSSR} 12:113 -- 122.

\bibitem[{Lee et~al.(2020)Lee, Paat, Stallknecht, \protect\BIBand{}
  Xu}]{LPSX2020}
Lee J, Paat J, Stallknecht I, Xu L (2020) Improving proximity bounds using
  sparsity. \emph{Proceedings of the 2020 International Symposium on
  Combinatorial Optimization} 115--127,
  \urlprefix\url{http://dx.doi.org/10.1007/978-3-030-53262-8_10}.

\bibitem[{Lee et~al.(2021)Lee, Paat, Stallknecht, \protect\BIBand{}
  Xu}]{LPSX2021}
Lee J, Paat J, Stallknecht I, Xu L (2021) Polynomial upper bounds on the number
  of differing columns of {$\Delta$}-modular integer programs.
  \emph{arXiv:2105.08160} .

\bibitem[{Lenstra(1983)}]{L1983}
Lenstra H (1983) Integer programming with a fixed number of variables.
  \emph{Mathematics of Operations Research} 8(4):538 -- 548.

\bibitem[{Mahler(1939)}]{M1939B}
Mahler K (1939) Ein {\"{u}}bertragungsprinzip f{\"{u}}r konvexe {K}{\"{o}}rper
  (in {G}erman). \emph{{\v{C}}asopis Pe{\v{s}}t. Mat. Fyz.} 68:93--102.

\bibitem[{N{\"a}gele et~al.(2022)N{\"a}gele, Santiago, \protect\BIBand{}
  Zenklusen}]{NSZ2021}
N{\"a}gele M, Santiago R, Zenklusen R (2022) Congruency-constrained {TU}
  problems beyond the bimodular case. \emph{In proceedings of SODA 2022}
  \urlprefix\url{http://dx.doi.org/10.1137/1.9781611977073.108}.

\bibitem[{Oertel et~al.(2020)Oertel, Paat, \protect\BIBand{}
  Weismantel}]{OPW2020}
Oertel T, Paat J, Weismantel R (2020) The distributions of functions related to
  parametric integer optimization. \emph{{SIAM} Journal on Applied Algebra and
  Geometry} 422--440, \urlprefix\url{http://dx.doi.org/10.1137/19M1275954}.

\bibitem[{Paat et~al.(2020)Paat, Weismantel, \protect\BIBand{}
  Weltge}]{PWW2020}
Paat J, Weismantel R, Weltge S (2020) Distances between optimal solutions of
  mixed-integer programs. \emph{Mathematical Programming} 179:455--468,
  \urlprefix\url{http://dx.doi.org/10.1007/s10107-018-1323-z}.

\bibitem[{Rudelson(2000)}]{R2000}
Rudelson M (2000) Distances between nonsymmetric convex bodies and the
  {MM}*-estimate. \emph{Positivity} 4:2:161--178.

\bibitem[{Santos(2012)}]{S2012}
Santos F (2012) A counterexample to the {H}irsch {C}onjecture. \emph{Annals of
  {M}athematics} 176:383--412,
  \urlprefix\url{http://dx.doi.org/10.4007/annals.2012.176.1.7}.

\bibitem[{Schrijver(1986)}]{AS1986}
Schrijver A (1986) \emph{Theory of linear and integer programming} (John Wiley
  \& Sons, Inc. New York, NY).

\bibitem[{Sturmfels(1996)}]{Sturmfels1995GrobnerBA}
Sturmfels B (1996) Gr{\"o}bner bases and convex polytopes (University Lecture
  Series, Volume 8, 162 pp.).

\bibitem[{Veselov \protect\BIBand{} Chirkov(2009)}]{VC2009}
Veselov S, Chirkov A (2009) Integer programming with bimodular matrix.
  \emph{Discrete Optimization} 6:220--222,
  \urlprefix\url{http://dx.doi.org/10.1016/j.disopt.2008.12.002}.

\bibitem[{Werman \protect\BIBand{} Magagnosc(1991)}]{W1991}
Werman M, Magagnosc D (1991) The relationship between integer and real
  solutions of constrained convex programming. \emph{Mathematical Programming}
  51:133--135, \urlprefix\url{http://dx.doi.org/10.1007/BF01586929}.

\bibitem[{Xu \protect\BIBand{} Lee(2019)}]{L2019}
Xu L, Lee J (2019) On proximity for $k$-regular mixed-integer linear
  optimization. \emph{Proceedings of WCGO 2019}, 438--447 (Springer),
  \urlprefix\url{http://dx.doi.org/10.1007/978-3-030-21803-4_44}.

\end{thebibliography}

\end{document}